\documentclass[11pt]{amsart}
\usepackage{latexsym}
\usepackage{amsfonts,amssymb,amsmath,amsthm,amscd}
\usepackage[mathscr]{eucal}

\newcommand{\su}{\mbox{${\mathfrak s \mathfrak u}$}}

\newcommand{\g}{\mbox{${\mathfrak g}$}}
\newcommand{\h}{\mbox{${\mathfrak h}$}}
\newcommand{\kf}{\mbox{${\mathfrak k}$}}

\newcommand{\p}{\mbox{${\mathfrak p}$}}

\newcommand{\tf}{\mbox{${\mathfrak t}$}}

\newcommand{\C}{\mbox{${\mathbb C}$}}

\newcommand{\I}{\mbox{${\mathbb I}$}}

\newcommand{\PP}{\mbox{${\mathbb P}$}}

\newcommand{\R}{\mbox{${\mathbb R}$}}
\newcommand{\Z}{\mbox{${\mathbb Z}$}}

\newcommand{\tr}{{\rm tr}}
\newcommand{\ric}{{\rm Ric}}

\newcommand{\wt}{${\rm conv}(\frac{1}{2}(d + {\mathcal W}))\,$}
\newcommand{\polyc}{$\Delta^{\bar{c}}\,$}
\newcommand{\polya}{$\Delta^{\bar{a}}\,$}
\newcommand{\ap}{\bar{a}^{\prime}}

\makeatletter
\def\numberwithin#1#2{\@ifundefined{c@#1}{\@nocnterrr}{%
  \@ifundefined{c@#2}{\@nocnterr}{%
  \@addtoreset{#1}{#2}%
  \toks@\expandafter\expandafter\expandafter{\csname the#1\endcsname}%
  \expandafter\xdef\csname the#1\endcsname
    {\expandafter\noexpand\csname the#2\endcsname
     .\the\toks@}}}}
\makeatother
\numberwithin{equation}{section}

\newtheorem{thm}[equation]{Theorem}

\newtheorem{prop}[equation]{Proposition}

\newtheorem{ex}[equation]{Example}

\newtheorem{rem}[equation]{Remark}
\newenvironment{rmk}{\begin{rem} \em}{\end{rem}}

\begin{document}

\title{Classifying Superpotentials: Three Summands Case}
\author{A. Dancer}
\address{Jesus College, Oxford University, Oxford, OX1 3DW, United
Kingdom}
\email{dancer@maths.ox.ac.uk}
\author{M. Wang}
\address{Department of Mathematics and Statistics,
McMaster University, Hamilton, Ontario, L8S 4K1, Canada}
\email{wang@mcmaster.ca}
\thanks{The second author is partially supported by NSERC Grant No. OPG0009421}

\date{revised \today}

\begin{abstract}
We give an overview of our classification results in \cite{DW4} and \cite{DW6}
for superpotentials of scalar curvature type of the cohomogeneity one Ricci-flat equations.
We then give an account of the classification in the case where the isotropy
representation of the principal orbit consists of exactly three distinct irreducible
real summands----the leftover case from \cite{DW6}.
\end{abstract}

\maketitle


\bigskip
\setcounter{section}{0}

\section{\bf The Classification Problem and Results}
This paper is a sequel to the papers \cite{DW4} and \cite{DW6}, in which
we studied the classification problem for superpotentials of the
cohomogeneity one Ricci-flat system. In this introductory section we will
give an overview of the problem and a summary of the results in the above
papers. We will then describe our new classification results for the case in
which the principal orbit contains exactly three irreducible summands---this
being the case that is left over from \cite{DW4} and \cite{DW6}.

Consider the Einstein equation $\ric(\bar{g}) = \Lambda \bar{g}$ for a Riemannian
metric $\bar{g}$ of cohomogeneity one with respect to some given compact Lie
group $G$ of isometries. We assume that the underlying $G$-manifold has
an interval as orbit space and that there is at least one singular orbit.
Over the part of the manifold consisting of principal orbits, we have
a family of equidistant $G$-homogeneous hypersurfaces. One can then
reduce the Einstein equations to a system of ordinary differential equations,
for example, by choosing a geodesic that intersects all principal
orbits orthogonally and using the arclength parameter along it as the independent
variable . The details of this reduction were first written down in \cite{BB}.
(See also \cite{EW} for information regarding smooth extension and local
existence near a singular orbit.) In the present paper, as in \cite{DW4} and
\cite{DW6}, we shall be concerned with the problem of finding interesting first
order subsystems of the Ricci-flat ODE system. One motivation for this problem
is that parallel, Killing, or holomorphic conditions are expressed by first
order systems.

One approach to finding interesting first order subsystems of the
Einstein ODEs  comes from viewing this system as a Hamiltonian
system with an extra constraint. In \cite{DW4} we gave a precise formulation
of such an approach. The Hamiltonian $\sf H$ is derived from the Einstein-Hilbert
lagrangian via a Legendre transformation, and is essentially the one used
in Hamiltonian formulations of General Relativity, cf \cite{Wa}, Appendix E.
It is interesting to note that even in the Riemannian case the kinetic energy part
of $\sf H$, which is the analogue of the Wheeler-de Wit metric, is also a
quadratic form of Lorentz signature. The potential energy part in turn depends on
the Einstein constant $\Lambda$ and the scalar curvature of the principal orbits.
The extra constraint is precisely the zero energy condition $\sf H = 0$.

For a general Hamiltonian system, we may define a {\it superpotential} as
a globally defined $C^2$ function $u$ on configuration space (with generalised
position variable $q$) such that
\begin{equation} \label{superpot}
   {\sf H}(q, du_q) = 0.
\end{equation}
In the physics literature the term superpotential has a more diffuse
meaning and hence different precise definitions in different situations.
One needs $u$ to be defined on all of configuration space because having
deduced the first order subsystem, one still has to find solutions
which extend smoothly to singular orbits and which define complete metrics.
Note also that the $C^2$ condition is imposed so that the map
$q \mapsto (q, du_q)$ is a $C^1$ graph and solutions of the
first order subsystem are solutions of the Einstein system.

In the context of the present paper, the configuration space is the cone
of $G$-invariant Riemannian metrics on the principal orbit, and its cotangent
space (with the canonical symplectic structure) is the momentum phase space.
We let $2r$ denote its dimension. By Proposition 1.1 of \cite{DW4}, for our
Einstein Hamiltonian system, the existence of a superpotential is equivalent to
the existence of $r$ functions on momentum phase space which depend linearly
on the conjugate momenta, whose common zero set lies in the constraint hypersurface,
and which Poisson commute upon restriction to the zero set. It follows that the
Hamiltonian vector field is tangent to this zero set, which is a lagrangian section
of the momentum phase space. The associated first order system is then the pull-back
of the Hamiltonian vector field to configuration space.

More explicitly, the first order system takes the form
\begin{equation} \label{FOS}
     \dot{q} = 2 v^{-1}(q) J^* \nabla u,
\end{equation}
where $v(q)$ is the volume of the metric $q$ on the principal orbit (relative
to some fixed background invariant metric), $J^*$ is the symmetric endomorphism
associated to the kinetic energy, and $\nabla$ is the Euclidean gradient on
configuration space. It is interesting to note that when the principal orbit
is not a torus, a superpotential cannot have any critical points
(cf \cite{DW4}, Proposition 1.3).

The attractiveness of the superpotential formalism lies in the uniform manner
in which the first order subsystem is derived. In recent years superpotentials
were used in string theory to derive the special holonomy conditions for non-compact
Ricci-flat Riemannian metrics with special asymptotics (see, for example, \cite{BGGG},
\cite{CGLP1}-\cite{CGLP3}). Roughly speaking, our classification theorems show that,
under appropriate assumptions that will be described more precisely below, besides
the first order systems which express these special holonomy conditions, there are
only a very small number of other possibilities. Included among the latter are the
first order subsystems arising from the doubly and triply warped products considered
in \cite{DW1} and \cite{DW3}. Thus while the detailed geometry associated to the
first order systems may change from case to case or may be absent, the existence
of the subsystem has a uniform characterization in terms of the
Hamiltonian/symplectic viewpoint.

From the partial differential equations point of view, the superpotential
equation (\ref{superpot}) is an implicitly defined first order equation that
can be interpreted as a time-independent Hamilton-Jacobi equation. Alternatively,
it can be transformed into an eikonal equation associated to a (variable)
Lorentz metric. In any case, from the general theory one certainly does not expect
the equation to have global regular solutions except in very special situations.
On the one hand, this is an indication that our classification problem would lead to
a finite number of possibilities. On the other hand, it also poses a challenge,
for there appears to be few clues in the literature about how to uncover these
nice situations.

Let $G/K$ be an $n$-dimensional connected homogeneous space where $G$ is a compact
Lie group and $K$ is a closed subgroup. Each such space may be viewed as the principal orbit of
a cohomogeneity one $G$-manifold $I \times (G/K)$ where $I$ is some interval in $\R$.
It determines a system of ODEs for Einstein metrics on $I \times (G/K)$ which we view
as a Hamiltonian system with constraint. The classification problem we pose is to find
all those $G/K$ whose Einstein system admits a superpotential. If possible one would also
like to determine all superpotentials and characterise the geometric conditions, if any,
they single out.

So far, we have examined this problem only in the Ricci-flat case, which is perhaps the
most interesting case. For a fixed $G$-invariant metric $g$ on $G/K$, its scalar curvature
${\sf S}_g$ is a constant. We will call the function $g \mapsto {\sf S}_g$ on configuration space
the {\it scalar curvature function} of $G/K$. In the case of a general $G/K$, this is
a complicated rational function (homogeneous of degree $-1$) of the components of $g$
(cf \cite{WZ}, \cite{BWZ}, \cite{Bo}). We will focus on the special case where
there are no multiplicities in the decomposition of the isotropy representation
of $G/K$ into irreducible real representations. In this case, the configuration space
is just $(\R_{+})^r$ and ${\sf S}_g$ is given by
(cf \cite{WZ}, Eq. (1.3))
\begin{equation} \label{SCF1}
     {\sf S} = \frac{1}{2} \sum_{i=1}^r \frac{a_i}{x_i}
     - \frac{1}{4} \sum_{i, j, k} \,[ijk] \frac{x_k}{x_i x_j},
\end{equation}
where $r$ is the number of irreducible summands, $x_i$ is the eigenvalue of $g$ (as
a symmetric automorphism with respect to some background normal homogeneous metric)
corresponding to the $i^{th}$ summand, and $[ijk] \geq 0$ is a coefficient depending
on the projection onto the $k^{th}$ summand of the Lie brackets of elements from
bases of the $i^{th}$ and $j^{th}$ summands. We let $d_i$ denote the (real) dimension
of the $i^{th}$ summand, so that $n= d_1 + \cdots + d_r$. Using exponential
coordinates defined by $x_i = e^{q_i}$ we may express $\sf S$ in the form
\begin{equation} \label{SCF2}
   {\sf S} = \sum_{w \in {\mathcal W}} A_w \, e^{w \cdot q},
\end{equation}
where $\mathcal W$ is a finite subset of $\Z^r \subset \R^r$ depending only on $G/K$
and $A_w$ are nonzero constants.

It follows from the above description that $\mathcal W$ consists of {\em weight
vectors} of the following three types:
\begin{enumerate}
\item[(i)] type I: one entry of $w$ is $-1$, the others are zero, with notation
         $(-1)^i$ where $i$ is the position of the non-zero entry,

\item[(ii)] type II: one entry is 1, two are -1, the rest are zero, with notation
          $(1^i, -1^j, -1^k)$ where $i, j$ and $k$ are the corresponding coordinate
          positions,

\item[(iii)] type III: one entry is 1, one is -2, the rest are zero, with analogous
       notation $(1^i, -2^j$).
\end{enumerate}

In this paper, we will call a function which is a finite linear combination
of exponentials in the $q_i$ a {\it function of scalar curvature type}. Here, as
well as in \cite{DW4} and \cite{DW6}, we will restrict ourselves to superpotentials
which are of scalar curvature type. While this is not a satisfactory assumption
from a general viewpoint, it does ensure that we are working with globally defined
superpotentials, and allows us to reduce the classification problem in this setting
to a problem involving convex polytopes as we shall explain below. Furthermore,
all known examples of superpotentials are of this type.

Let us then write the superpotential $u$ in the form
\begin{equation} \label{potform}
     u = \sum_{\bar{c} \in {\mathcal C}} F_{\bar{c}} \,\, e^{{\bar{c}}\cdot q}
\end{equation}
where $\mathcal C$ is a finite subset of $\R^r$ to be determined, and
$F_{\bar{c}}$ are nonzero unknown real constants. By \cite{DW2}, the Hamiltonian is given
by
\begin{equation} \label{ham}
   {\sf H} = v^{-1}(q) J(p,p) + v(q)((n-1) \Lambda - {\sf S}_q)
\end{equation}
where $v$ and ${\sf S}_q$ are respectively the relative volume and
scalar curvature function mentioned above, $\Lambda$ is the Einstein constant,
and $J(p,p)$ is the
quadratic form
\begin{equation} \label{Jform}
   J(p,p) = \frac{1}{n-1} \left(\sum_{i=1}^r p_i \right)^2 - \sum_{i=1}^r \frac{p_i^2}{d_i}
\end{equation}
with signature $(1, r-1)$.

Substituting (\ref{potform}) into the superpotential equation (with $\Lambda = 0$)
we obtain, for each $\xi \in \R^r$, the equation
\begin{equation} \label{eqnF}
\sum_{\bar{a}+\bar{c} = \xi} J(\bar{a},\bar{c}) \ F_{\bar{a}} F_{\bar{c}} = \left\{
  \begin{array}{ll}
    A_w  &   \mbox{if} \ \xi = d + w \  \mbox{for some} \ w \in \mathcal{W} \\
     0    &  \mbox{if} \ \xi \notin d + \mathcal W
   \end{array} \right.
\end{equation}
where $d$ denotes the vector $(d_1, \cdots,d_r)$. Hence the classification
problem requires us to find all $G/K$ (satisfying our hypotheses) for which
this set of equations has a solution, and in each such case to find all
the solutions.

In view of the above formulation of the problem, if the weight vector $w \in
{\mathcal W}$, then there are elements $\bar{a}, \bar{c} \in {\mathcal C}$
such that $d+w = \bar{a} + \bar{c}$, and so $\frac{1}{2}(d+w)$ lies in the
convex hull of $\mathcal C$. Hence the pair of convex hulls
\begin{equation} \label{incl}
 {\rm conv}(\frac{1}{2}(d + {\mathcal W})) \subset {\rm conv}(\mathcal C)
\end{equation}
plays a critical role in our analysis. $\mathcal C$ is of course unknown, but
the possibilities for $\mathcal W$, which depends on $G/K$, are also to be determined.
We will now make a further assumption, namely that ${\rm conv}(\mathcal W)$
has dimension $r-1$. As was shown in the proof of Theorem 3.11 in \cite{DW2},
this assumption is satisfied if $G$ is semisimple.

The case $r=1$ is very special. In this situation, $G/K$ is isotropy irreducible,
so the cohomogeneity one space cannot have a singular orbit. As well, the quadratic
form $J$ is positive definite. In any event, there is always a superpotential
of the desired type, as was shown at the end of Section 1 of \cite{DW4}. Henceforth
we shall assume that $r \geq 2$. Unless otherwise stated, we will use the
Lorentz metric $J$ on $\R^r$.

There are now two main cases, depending on whether or not the inclusion in
(\ref{incl}) is strict. As we remarked in Section 1 of \cite{DW6}, this condition
is equivalent to the existence of a null vertex in $\mathcal C$.
Indeed, (\ref{eqnF}) implies that a non-null vertex $\bar{c}$ of ${\rm conv}(\mathcal C)$
must be of the form $d+w$ for some $w \in {\mathcal W}$. On the other hand,
if $\bar{c}$ is a null vertex and equality in (\ref{incl}) holds, then $2\bar{c} = d+ w$
for some $w \in {\mathcal W}$, and (\ref{eqnF}) then fails for $\xi = d + w$.

The case of equality is the less complicated of the two cases, for we need only
to analyse (\ref{eqnF}) in terms of the convex polytope \wt.
This was done in \cite{DW4} where we proved

\begin{thm} \label{nonnullthm}
Let $G$ be a compact connected Lie group and $K$ a closed connected subgroup
such that the isotropy representation of $G/K$ decomposes into a sum of $r$
pairwise inequivalent irreducible real summands. Assume that ${\rm conv}(\mathcal W)$
has dimension $r-1$.

Suppose that the cohomogeneity one Ricci-flat equations with principal orbit
$G/K$ admit a superpotential of scalar curvature type $($\ref{potform}$)$ such that all
elements of $\mathcal C$ are non-null. Then the possibilities, up to permutations of
the irreducible summands, are given by
\begin{enumerate}
\item ${\mathcal W} = \{ (-1) \}$ and $G/K$ is isotropy irreducible,

\item ${\mathcal W} = \{(1, -2), (0, -1)\},$  and
       $G/K = (SO(3) \times SO(2))/\Delta SO(2),$

\item $\mathcal{W} = \{ (1, -2), (-1,0), (0, -1) \},$
       and $G/K$ is one of $(SU(3) \times SO(3))/ \Delta SO(3) \approx SU(3)$,
       $(Sp(2) \times Sp(1))/(Sp(1) \times \Delta Sp(1)) \approx S^7$,
$(SU(3) \times SU(2))/(U(1) \cdot \Delta SU(2))\approx SU(3)/U(1)_{11},$
 or $Sp(2)/(U(1)\times Sp(1))= SO(5)/U(2)\approx \C \PP^3, $

\item $\mathcal{W} = \{ (1,-2,0), (1,0,-2), (0,1,-2), (0,-1,0), (0,0,-1) \},$
     and $G/K$ is $S^7$ written in the form
    $(Sp(2) \times U(1))/(Sp(1)\Delta U(1))$,

\item $\mathcal{W}= \{(1,-1,-1), (-1,1,-1), (-1,-1,1), (-1,0,0), (0,-1,0), (0,0,-1)\}$,
      $d=(2,2,2)$  and $G/K = SU(3)/T$,

\item $ \mathcal{W} = \{ (1,-2,0,0), (1,0,-2,0),
 (1,0,0,-2),(0,1,-1,-1), (0,-1,1,-1),(0,-1,-1,1),$ \\
 $ (0,-1,0,0), (0,0,-1,0), (0,0,0,-1) \},$ and $G/K$
 is an Aloff-Wallach space $SU(3)/U_{kl}$, where
 $U_{kl}$ denotes the circle subgroup consisting of the diagonal
 matrices ${\rm diag}(e^{ik\theta}, e^{il\theta}, e^{im\theta})$
 with $k+l+m = 0, (k,l)=1$, and $\{k, l, m\} \neq \{1, 1, -2\}$
or $\{1, -1, 0\}$,

\item a local product of an example in $($1$)$ $(n>1)$, $($3$)$, or $($5$)$
     with a circle.
\end{enumerate}
In all of the above cases, there is a superpotential of scalar curvature
type that is unique up to an overall minus sign and an additive constant.
\end{thm}

\begin{rmk} \label{rmk-nonnull}

(a) In the above theorem, the first order subsystem resulting from case (2)
corresponds to the hyperk\"ahler condition. The first order subsystems
of the second and fourth subcases in (3) as well as those of (4), (5), (6)
correspond to special holonomy $G_2$ or ${\rm Spin}(7)$ according to
whether $n=6$ or $7$.  The first and third subcases of (3) do not allow the
addition of a singular orbit and are not related to special holonomy.

(b) We refer the interested reader to Section 7 of \cite{DW4} for further
information about the superpotentials and solutions of the first order systems.

(c) Case (7) results from a general property of superpotentials
of scalar curvature type without null weight vectors associated with a
principal orbit $G/K$ having no trivial summands in its isotropy representation
(cf Remark 2.8 in \cite{DW4}).
\end{rmk}

In the situation where $\mathcal C$ contains a null vector, our analysis
is more complicated because we have to study the relative positions of
${\rm conv}(\mathcal C)$ and ${\rm conv}(\frac{1}{2}(d+ {\mathcal W}))$.
To do this, we take a generic hyperplane separating a null vertex
$\bar{c}$ from ${\rm conv}(\frac{1}{2}(d+ {\mathcal W}))$ and examine the
image $\Delta^{\bar{c}}$ of the perspective projection of \wt onto the
separating hyperplane. Analysing the vertices and edges of $\Delta^{\bar{c}}$
is equivalent to analysing the edges and $2$-dimensional faces of
${\rm conv}(\frac{1}{2}(d+ {\mathcal W}))$. The large number of possibilities
of the $2$-dimensional faces makes our task rather onerous (cf \cite{DW5} for
the 246+ faces corresponding just to the edges in $\Delta^{\bar{c}}$ connecting what
we refer to in \cite{DW6} as adjacent type (1B) vertices). However, if $r \geq 4$, i.e.,
when \wt is at least $3$-dimensional, it
turns out that we can avoid analysing in detail all the convex subshapes of the
particular hexagonal face we denoted by (H1) in \cite{DW6}. The main result
of that paper is

\begin{thm} \label{nullthm}
Let $G$ be a compact connected Lie group and $K$ a closed connected
subgroup such that the isotropy representation of $G/K$ is the direct
sum of $r$ pairwise inequivalent $\R$-irreducible summands.
Assume that  $\dim {\rm conv}(\mathcal W) = r-1$.

Suppose the cohomogeneity one Ricci-flat equations with $G/K$ as
principal orbit admit a superpotential of scalar curvature type $($\ref{potform}$)$
where $\mathcal C$ contains a null vertex. Then either $r \leq 3$, or, up to permutations
of the irreducible summands, we have
\begin{eqnarray*}
 {\mathcal W} &=& \{ (-1)^i, (1^1, -2^i): 2 \leq i \leq r \}, \ \
   d_1=1, \\
 {\mathcal C} &=& \frac{1}{2}(d+ \{(-1^1), (1^1, -2^i): 2 \leq i \leq r \})
  \ \  \mbox{\rm with} \ \  r \geq 2,
\end{eqnarray*}
and the superpotential of scalar curvature type is unique up to an overall minus sign
and an additive constant.
\end{thm}

\begin{rmk}
Note that the possibility described in the above theorem is realised by circle bundles
over a product of $r-1$ Fano (homogeneous) K\"ahler-Einstein manifolds
The corresponding first order subsystem corresponds to the Calabi-Yau condition.
More discussion of this example can be found in Section 8 of \cite{DW4}.
\end{rmk}

The $r=2$ case was treated at the end of \cite{DW6}. See Proposition \ref{c-points}
and the remarks following it for further details. This paper will be devoted to
treating the remaining case $r=3$. We summarise the $r=2,3$ results in the following

\begin{thm} \label{smallrthm}
Let $G$ be a compact Lie group and $K$ be a closed subgroup such that
$G/K$ is connected and $\kf$ is not a maximal ${\rm Ad}_K$-invariant subalgebra
of $\g$. Assume that the isotropy representation of $G/K$ splits into $2$ or $3$
pairwise distinct irreducible real summands and $\mathcal C$ contains a null vector.

If the Ricci-flat cohomogeneity one Einstein equations with $G/K$ as principal
orbit admit a superpotential of scalar curvature type, then, up to permutations of the
irreducible summands, the possibilities are
\begin{enumerate}
\item ${\mathcal W} = \{(-1, 0, 0), (0, -1, 0), (0, 0, -1)\}$ with
       $ d=(3,3,3), (2, 4, 4),$ or $(2, 3, 6)$.
\item ${\mathcal W} = \{(0, -1, 0), (0, 0, -1), (1, -2, 0), (1, 0, -2) \}$ with
      $d_1 =1$,
\item ${\mathcal W } = \{(-1, 0), (0, -1)\}$, with $\frac{4}{d_1} + \frac{1}{d_2} = 1$,
\item ${\mathcal W } = \{(0, -1), (1, -2)\}$, with either $d_1 = 1$ or $\frac{4}{d_1}+
     \frac{9}{d_2} = 1$.
\end{enumerate}
In each of the above cases, there is a superpotential of scalar curvature type that is
unique up to a sign and an additive constant.
\end{thm}

\begin{rmk} \label{rmksmallr}
(a) Case (2) and the $d_1=1$ subcase of Case (4) in the above theorem are
respectively just the $r=3, 2$ cases of the Calabi-Yau case in Theorem \ref{nullthm}.
Furthermore, the second possibility is realised by the complete, non-compact
B{\'erard} Bergery examples \cite{BB}.

(b) Case(3) and the second subcase of Case (4) are realised by the {\em explicit}
doubly-warped examples studied in \cite{DW1}. Case (1) is realised by the
triply-warped examples studied in \cite{DW3} and \cite{DW4}. The first order
subsystems for the $d=(3,3,3)$ and $(2,4,4)$ are integrable by quadratures. Further
details can be found in Section 8 of \cite{DW4}.
\end{rmk}

Note that in Theorems (\ref{nonnullthm}) and (\ref{nullthm}) we have assumed the
connectedness of both $G$ and $K$ instead of the connectedness of $G/K$ as in
Theorem (\ref{smallrthm}). If we drop the more stringent condition, further examples
of superpotentials arise. Noteworthy examples without null weights, described
 in more detail in Section 7 of \cite{DW4}, include
\begin{enumerate} \label{disconnectedcases}
\item  $G/K = O(3)/(O(1)\times O(1) \times O(1))$, for which there are two
   {\em distinct} superpotentials of scalar curvature type,
\item  $G/K = ([SU(2)\times SU(2)\times \Delta U(1)] \ltimes \Z_2 )/ (\Delta U(1) \times \Z_2)
     \approx S^3 \times S^3,$ for which a superpotential was found in \cite{BGGG}
     and \cite{CGLP2} in connection with $G_2$-holonomy.
\end{enumerate}
On the other hand, large portions of \cite{DW6} and (to a lesser extent) \cite{DW4}
do not require $G$ and $K$ to be connected, as we have indicated at appropriate
points of those papers. It would be interesting to remove this condition
in the classifications, as the above examples show. Therefore, Theorem \ref{smallrthm}
can be viewed as a contribution to that effort, as well as the completion of
Theorem \ref{nullthm}.

\section{\bf Overview of the three-summands case}

In this section we first recall some generalities about the classification
for the null case described in \cite{DW6} and adapt them to the $r=3$ case.
We then give an overview as to how the classification will proceed.

Recall that $G$ and $K$ are compact Lie groups such that $G/K$ is connected
and almost effective. We choose a bi-invariant metric on $G$ and take the
induced normal metric on $G/K$ as our background metric. Let
$$ \g = \kf \oplus \p_1 \oplus \p_2 \oplus \p_3 $$
be the associated orthogonal decomposition where $\p_i$ are the irreducible
real summands of the isotropy representation. We assume that these ${\rm Ad}_K$
modules are pairwise inequivalent and that $\kf$ is not a maximal ${\rm Ad}_K$
invariant subalgebra in $\g$.

The following are some useful facts about the scalar curvature function of
$G/K$:
\begin{enumerate}
\item[(a)] For a type I vector $w$, the coefficient $A_w > 0$
   in (\ref{SCF2}) while for type II and type III vectors, $A_w < 0$.

\item[(b)] The type I vector with $-1$ in the $i$th position is absent from $\mathcal W$
iff the corresponding summand $\p_i$ is an abelian subalgebra which
satisfies $[\kf, \p_i] = 0$ and $[\p_i, \p_j] \subset \p_j$ for
all $j \neq i$. If the isotropy group $K$  is connected, these
last conditions imply that $\p_i$ is $1$-dimensional, and the
$\p_j, j \neq i,$ are irreducible representations of the (compact)
analytic group whose Lie algebra is $\kf \oplus \p_i$.

\item[(c)] If $(1^i, -1^j, -1^k)$ occurs in $\mathcal W$ then its permutations
$(-1^i, 1^j, -1^k)$ and $(-1^i, -1^j, 1^k)$ do also.

\item[(d)] If $\dim \p_i = 1$ then no type III vector with $-2$ in place $i$ is
 present in $\mathcal W$.

\item[(e)] If $I$ is a subset of $\{1, \cdots, r \},$ then
each of the equations $\sum_{i \in I} x_i = 1$ and $\sum_{i \in I} x_i =-2$
defines a face (possibly empty) of conv$(\mathcal W)$.
In particular, all type III vectors in $\mathcal W$ are vertices
and $(-1^i, -1^k, 1^j) \in {\mathcal W}$ is a vertex unless
both $(-2^i, 1^j)$ and $(-2^k, 1^j)$ lie in $\mathcal W$.

\item[(f)] For $v, w \in \mathcal{W}$ (or indeed for any $v,w$ such that $\sum v_i$
or $\sum w_i=-1$), we have
\begin{equation} \label{Jexpr}
J(v +d, w+ d) = 1 -\sum_{i=1}^{r} \frac{v_i w_i}{d_i}.
\end{equation}
\end{enumerate}

We have the following schematic picture of the full set of weight vectors
which can appear in the scalar curvature function (\ref{SCF2}) when $r=3$.
This is the hexagon (H1) mentioned in Section 6 of \cite{DW6}.

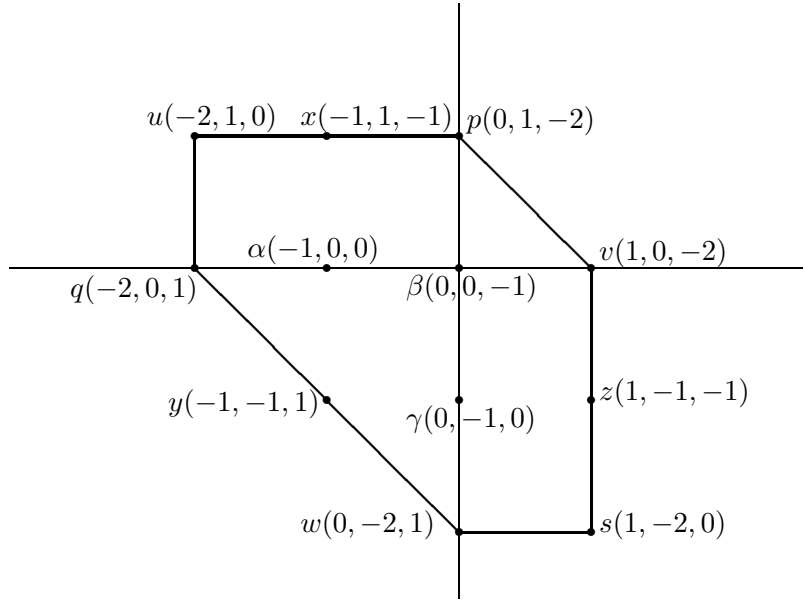
\begin{picture}(300, 200)(0, 10)
\put(150, 100){\circle*{3}}
\put(200, 100){\circle*{3}}
\put(250, 100){\circle*{3}}
\put(300, 100){\circle*{3}}
\put(150, 150){\circle*{3}}
\put(200, 150){\circle*{3}}
\put(250, 150){\circle*{3}}
\put(200, 50){\circle*{3}}
\put(250, 50){\circle*{3}}
\put(300, 50){\circle*{3}}
\put(250, 0){\circle*{3}}
\put(300, 0){\circle*{3}}
\put(80, 100){\line(1, 0){300}}
\put(250, -25){\line(0, 1){225}}
\thicklines
\put(150, 150){\line(1, 0){100}}
\put(150, 100){\line(0, 1){50}}
\put(300, 100){\line(-1,1){50}}
\put(300, 0){\line(0,1){100}}
\put(300, 0){\line(-1,0){50}}
\put(250, 0){\line(-1, 1){100}}
\put(303, 103){$v(1, 0, -2)$}
\put(303, 50){$z(1, -1, -1)$}
\put(303, 0){$s(1, -2, 0)$}
\put(253, 153){$p(0, 1, -2)$}
\put(190, 154){$x(-1, 1, -1)$}
\put(132, 154){$u(-2, 1, 0)$}
\put(103, 89){$q(-2, 0, 1)$}
\put(140, 46){$y(-1, -1, 1)$}
\put(190, 0){$w(0, -2, 1)$}
\put(170, 104){$\alpha(-1, 0, 0)$}
\put(230, 90){$\beta(0, 0, -1)$}
\put(230, 40){$\gamma(0, -1, 0)$}
\end{picture}

\vspace{2cm}

Recall that we are assuming that ${\rm conv}(\mathcal W)$ has dimension
$r-1=2$, and so the scalar curvature function determines a $2$-dimensional
convex subpolygon of the above hexagon. By Theorem 1.5 in \cite{DW6} we know
that the set of weights of a superpotential (of scalar curvature type) lies in
the (hyper)plane $\{\bar{x}: \sum_i \bar{x}_i = \frac{1}{2}(n-1) \}$, which contains
\wt  as well. Unless otherwise stated, we shall
be working in this plane. We will denote vectors in this plane by $\bar{c}, \bar{u},
\ldots$, and the corresponding vectors in the plane $x_1 + x_2 + x_3 = -1$ by
$c, u, \ldots.$

Henceforth we shall assume that we are not in the Calabi-Yau case, i.e.,
case (2) of Theorem \ref{smallrthm}. Recall that by Theorem 3.14 in \cite{DW6}
we are in this case if $\mathcal C$ contains a type I vector.

Let $\bar{c}$ denote a null vector in $\mathcal C$. As in Section 3 of \cite{DW6},
we choose a generic affine line separating $\bar{c}$ from \wt. The image of
\wt under the perspective projection through $\bar{c}$ is a segment $\Delta^{\bar{c}}$
in the affine line. The endpoints of this segment are of three types
(cf Theorem 3.8 of \cite{DW6}):
\begin{enumerate}
\item[(i)] Type (1A)-- the vertex is orthogonal to $\bar{c}$,
\item[(ii)] Type (1B)-- the line connecting $\bar{c}$ to the vertex intersects
   \wt in a unique point $\bar{x}$ and there is a null vector $\bar{a} \in {\mathcal C}$
   such that $2\bar{x} = \bar{a} + \bar{c}$,
\item[(iii)] Type (2)-- the vertex is not orthogonal to $\bar{c}$ and the line
joining it and $\bar{c}$ intersects \wt in an edge $\bar{v}\bar{w}$.
\end{enumerate}

The endpoints cannot both be of type (1A) since $\bar{c}^{\perp} \cap$ \wt
has dimension $\leq 1$. They also cannot both be of type (2) by Theorem 7.1
of \cite{DW6}, the fact that $r=3$, and the assumption that we are not
in the Calabi-Yau case. If one endpoint is type (2) and the other is type (1A),
then since $r=3$, the results in Section 9 of \cite{DW6} show first that
there must be points of $\frac{1}{2}(d+ {\mathcal W})$ in the interior of
the edge $\bar{v}\bar{w}$ mentioned in (iii) in the previous paragraph.
Next the parts of the proof of Theorem 9.2 in \cite{DW6} dealing with this
situation show that there cannot be an endpoint of type (1A), a contradiction.

Therefore, one of the endpoints of $\Delta^{\bar{c}}$ must be of type (1B).
If the other endpoint is of type (2), we shall show in Section 4 that $c$
must be a type I weight vector, so by Theorem 3.14 in \cite{DW6} we are in the
Calabi-Yau case. In Section 5 we will show that if both endpoints of $\Delta^{\bar{c}}$
are of type (1B) then we can only be in the situation of case (1) of
Theorem \ref{smallrthm}. The analysis required for this section is the most
onerous of all. Finally, we shall show in Section 6 that the other endpoint
of $\Delta^{\bar{c}}$ also cannot be of type (1A), and this then completes
the proof Theorem \ref{smallrthm}.

\section{\bf Listing of Subpolygons}

In this section we compile the list of convex subpolygons of the hexagon
(H1) which may occur as ${\rm conv}(\mathcal W)$. The properties (a)-(f)
in the previous section for a scalar curvature function will be used
without mention below. We will also refer freely to the labelled points in (H1).

Since we are assuming that $\kf$ is not a maximal ${\rm Ad}_K$ invariant
subalgebra in $\g$, after permuting the irreducible summands $\p_i$, we have
the following possibilities:

\noindent{I. \,\, $\kf \subset \h_1 \subset \h_2 \subset \g$} where $\h_1 = \kf \oplus \p_1$,
$\h_2 = \h_1 \oplus \p_2$ are ${\rm Ad}_K$ invariant intermediate subalgebras.

Since $[\p_1,\p_1] \subset \kf \oplus \p_1$, it follows that $(-2, 1, 0)$
and $(-2, 0, 1) \notin {\mathcal W}$. Next, $[\p_1, \p_2] \subset \p_2$
implies that $\mathcal W$ does not contain any type II vectors. Finally,
$[\p_2, \p_2] \subset \kf \oplus \p_1 \oplus \p_2$ means that $(0, -2, 1) \notin
{\mathcal W}$. We therefore obtain the trapezoid $pvs\alpha$ (Q6 below) with interior point
$\beta$ and the midtpoint $\gamma$ of $s\alpha$. This trapezoid does occur
as ${\rm conv}(\mathcal W)$ of some $G/K$, see Example 7 of \cite{WZ}

\noindent{II. \, \, $\kf \subset \h \subset \g$ where $\h = \kf \oplus \p_1 \oplus \p_2$}
but neither $\kf \oplus \p_1$ nor $\kf \oplus \p_2$ is an ${\rm Ad}_K$ invariant subalgebra.

Since $[\p_i, \p_3] \subset \p_3$ for $i=1, 2$ there are no type II vectors.
Also, $\h$ being a subalgebra means that $(-2, 0, 1), (0, -2, 1) \notin {\mathcal W}$.
We now obtain the trapezoid $pvsu$ (Q4 below) together with all the points $\alpha, \beta,$
and $\gamma$.

\noindent{III. \, \, $\kf \subset \h \subset \g$} where $\h = \kf \oplus \p_1$
but $\h \oplus \p_i, i=2, 3,$ do not form an ${\rm Ad}_K$-invariant subalgebra.

Since $[\p_1, \p_1] \subset \h$, $(-2, 1, 0), (-2, 0, 1) \notin {\mathcal W}$.
So the possible polygons are (up to permutations of the summands) the basic hexagon $pvswyx$
(H2 below) and the subpolygons $pvzwyx$ (H3 below), $pvsw\alpha$ (P1 below) and $pvsyx$
(P2 below).

We must next consider all admissible subpolygons of the above polygons. After taking
into account permutations of the irreducible summands, we arrive at the following list:

\smallskip

\noindent{\bf Triangles:}
\begin{enumerate}
\item[(T1)] vertices $\alpha, \beta, \gamma$ (case of triply warped products)
\item[(T2)] vertices $p, v, \beta$
\item[(T3)] vertices $\beta, \gamma, v$
\item[(T4)] vertices $p, v, \gamma$ with extra point $\beta$
\item[(T5)] vertices $p, s, \gamma$ with extra point $\beta$
\item[(T6)] vertices $v, s, \alpha$ with extra points $\beta, \gamma$
\item[(T7)] vertices $\alpha, \beta, s$ with point $\gamma$
\item[(T8)] vertices $x, y, z$ and all the midpoints of the edges (realised by ${\rm SU}(3)/T$)
\item[(T9)] vertices $p, v, w$ with points $\beta, \gamma$
\item[(T10)] vertices $p, w, \alpha$ with points $\beta, \gamma$
\item[(T11)] vertices $v, s, u$ with interior point $\beta$ and points $\alpha, \gamma$
\item[(T12)] vertices $\alpha, p, s$ with interior point $\beta$ and point $\gamma$
\end{enumerate}

\noindent{\bf Quadrilaterals:}
\begin{enumerate}
\item[(Q1)] parallelogram with vertices $p, v, s, \gamma$ and point $\beta$
\item[(Q2)] square with vertices $p, v, \gamma, \alpha$ with interior point $\beta$
\item[(Q3)] rectangle with vertices $v, s, y, x$, interior points $\beta, \gamma$ and
       points $\alpha, z$
\item[(Q4)] trapezoid with vertices $p, v, s, u$, interior point $\beta$ and points $\alpha, \gamma$
\item[(Q5)] trapezoid with vertices $v, s, w, p$ and points $\beta, \gamma$
\item[(Q6)] trapezoid with vertices $p, v, s, \alpha$ with interior point $\beta$ and point $\gamma$
\item[(Q7)] trapezoid with vertices $v, s, \gamma, \beta$ (Calabi-Yau case)
\item[(Q8)] trapezoid with vertices $p, z, y, x$ with interior point $\beta$ and points $\alpha, \gamma$
\item[(Q9)] irregular quadrilateral with vertices $p, v, w, \alpha$ and interior points $\beta, \gamma$
\end{enumerate}

\noindent{\bf Pentagons:}
\begin{enumerate}
\item[(P1)] vertices $p, v, s, w, \alpha$ with interior points $\beta, \gamma$
\item[(P2)] vertices $p, v, s, y, x$ with interior points $\beta, \gamma$ and points $\alpha, z$
\item[(P3)] vertices $p, v, z, y, x$ with interior point $\beta$ and points $\alpha, \gamma$
\item[(P4)] vertices $p, z, s, y, x$ with interior points $\beta, \gamma$ and point $\alpha$
\item[(P5)] vertices $p, z, w, y, x$ with interior points $\beta, \gamma$ and point $\alpha$
\end{enumerate}

\noindent{\bf Hexagons:}
\begin{enumerate}
\item[(H2)] vertices $p, v, s, w, y, x$ with interior points $\beta, \gamma$ and points $\alpha, z$
\item[(H3)] vertices $p, v, z, w, y, x$ with interior points $\beta, \gamma$ and point $\alpha$
\end{enumerate}

Note that the polygons which do not contain all the points $\alpha, \beta, \gamma$ are
(T2)-(T5), (T9), (Q1), (Q7), and (Q5). In the analysis for null vertices $\bar{c} \in  {\mathcal C}$
for which the endpoints of $\Delta^{\bar{c}}$ are both of type (1B), these cases were
already considered in \cite{DW6} (cf comments just before Theorem 6.12).

\section{\bf Case of Endpoints of Type (1B)(2)}

Let $\bar{c} =\frac{1}{2}(d+c)$ be a null vertex in $\mathcal C$ such that
the endpoints of $\Delta^{\bar{c}}$ are of type (1B) and (2). We will show
in this section that this possibility cannot occur.  Let $\bar{V}\bar{W}$
denote the edge of \wt collinear with the ray from $\bar{c}$ to the
type (2) endpoint, with $\bar{V}$ denoting the point closer to $\bar{c}$.
By Theorem 9.2 in \cite{DW6} there must be points of $\frac{1}{2}(d+{\mathcal W})$
in the interior of $\bar{V}\bar{W}$. The possible configurations are given
in Table 3 in \cite{DW6}, where $v, w$ there are now denoted by $V, W$
to avoid conflict with the labels in (H1).

We now consider the polygons listed in the previous section having
edges which contain interior points lying in $\frac{1}{2}(d+{\mathcal W})$.
Note that this eliminates the polygons (T1)-(T3), (Q2), (Q9) and (P1).
Besides the nullity of $\bar{c}$, further conditions result by considering
the (1B) endpoint. Recall also that if \wt lies in one of the open half-planes
of an affine line, then there is at most one element of $\mathcal C$, necessarily
null, lying in the opposite half-plane (cf Proposition 3.3 of \cite{DW6}).
For a polygon which does not contain all three type I vectors,
we note that (cf paragraph before Theorem 6.12 in \cite{DW6})
$\mathcal C$ cannot contain a null element $\bar{a}$ whose
associated segment $\Delta^{\bar{a}}$ has both endpoints of type (1B).

We shall also need the following extension of Theorem 3.11 in \cite{DW6}.

\begin{prop} \label{c-points}
Let $\bar{c} \in {\mathcal C}$ be a null vector such that \polyc  \,
has a type $($2$)$ vertex. Let $\bar{V}\bar{W}$ denote the associated
edge in \wt $($with $\bar{V}$ nearer to $\bar{c}$$)$. Assume that
the midpoint $\alpha$ of $\bar{V}\bar{W}$ lies in $\frac{1}{2}(d+{\mathcal W})$.
Then $\bar{W} \notin {\mathcal C}$, i.e., there must be an element of
$\mathcal C$ lying beyond $\bar{W}$.

If we write $\bar{c}= \lambda \bar{V} + (1-\lambda)\bar{W}$, we have
$1<\lambda \leq \frac{3}{2}$. Furthermore, we cannot be in case $($2$)$,
$($3$)$, or $($4$)$ with $d_1 \neq d_3$ in Table 3 of \cite{DW6}.
\end{prop}

\begin{proof} Let $\bar{c}=\bar{c}^{(0)},\cdots, \bar{c}^{(m+1)}, m \geq 1,$
denote the points lying in $\mathcal C \cap \bar{V}\bar{W}$, ordered consecutively.
By Remark 3.12 in \cite{DW6}, we have $2\bar{V}=\bar{c}^{(0)} + \bar{c}^{(1)}$.
If for some $j>1,$ $\bar{c}^{(0)} + \bar{c}^{(j)} \in d+{\mathcal W}$, it follows
that $\bar{c}^{(0)} + \bar{c}^{(j)}= 2\bar{\alpha}$ or $2\bar{W},$ and $j=m+1$
since $\bar{c}^{(j)}$ would have to be null. Indeed the second possibility cannot
hold since by symmetry we also must have $\bar{c}^{(m)} + \bar{c}^{(m+1)} = 2\bar{W},$
which implies that $V=W$.

If $\bar{W} \in {\mathcal C}$,  then for all $j>1$, $\bar{c}^{(0)} +
\bar{c}^{(j)} \notin d+{\mathcal W}$. Furthermore, since now $\bar{W}=\bar{c}^{(m+1)}$
we have $J(\bar{c}^{(m)}, \bar{W})=0$ (otherwise $\bar{c}^{(m)}+\bar{c}^{(m+1)} = 2\bar{v}
$ or $2\bar{\alpha}$, neither of which are possible). By arguments very similar to
those in part (C) of the proof of Theorem 3.11 in \cite{DW6}, we are reduced to
the cases $m=0$ or $m=1$. However, in either case, the term in the scalar curvature
function corresponding to $\alpha$ is unaccounted for.

Hence $\bar{c}^{(m+1)}$ is null and $\bar{c}^{(0)} + \bar{c}^{(m+1)} = 2\bar{\alpha}$
holds; otherwise for all $j>1$ we have $\bar{c}^{(0)}+\bar{c}^{(j)} \notin d+{\mathcal W}$
and arguments very similar to those in part (B) of the proof of Theorem 3.11 in \cite{DW5}
lead to a contradiction. Since $\alpha$ is the midpoint of $VW$, this last condition
implies that $\bar{\alpha}$ is the midpoint of $\bar{c}^{(1)}\bar{c}^{(m)}$.
This yields the estimate $\lambda \leq 3/2$.

For the remaining assertions of the theorem, we look at Table 3 of \cite{DW6}.
Note that case (1) of that table does not satisfy the midpoint hypothesis of our
theorem. Let us consider $(1-t)V + tW$ where $V$ (resp. $W$) is the same as $v$
(resp. $w$) in the table. Let $t_1, t_2$ be the values corresponding
to $\bar{c}^{(0)}$ and $\bar{c}^{(m+1)}$ respectively. For case (3), we obtain
the condition $t_1 +t_2 = \frac{d_2}{d_1+d_2}$ from the null conditions for
$\bar{c}^{(0)}$ and $\bar{c}^{(m+1)}$. From the condition
$\bar{c}^{(0)} + \bar{c}^{(m+1)} = 2\bar{\alpha},$ we get instead
$t_1+t_2 =1$. So case (3) (and hence case (2)) cannot occur.
For case (4) we obtain $t_1 +t_2 =1$, which is equivalent to $d_1=d_3$.
So if this does not hold, case (4) cannot occur. This completes the proof
of the proposition.
\end{proof}

\begin{rmk} \label{2summands}
The above proposition has an immediate application to the null case of the
$r=2$ classification since in that case \polyc reduces to a single point and
must be of type (2). As we are assuming that $\kf$ is not a maximal ${\rm Ad}_K$
invariant subalgebra in $\g$, case (1) of Table 3 in \cite{DW6} does not
occur. This leaves cases (2),(3) since $r=2$, but these are eliminated by
the above proposition. Therefore we cannot have interior points of $\bar{V}\bar{W}$
lying in $\frac{1}{2}(d+ {\mathcal W})$ and we are in the last two possibilities of
Theorem \ref{smallrthm} by Theorem 9.2(i) of \cite{DW6}.
\end{rmk}

Using Proposition \ref{c-points} we now examine each of the polygons in Section 3
for a null $\bar{c} \in {\mathcal C}$ with a type (1B) and a type (2) endpoint
in \polyc. In fact we can immediately eliminate (T4)-(T7), (T12), (Q1) and (Q6).
Another useful fact is that since there must be null elements of $\mathcal C$
beyond either endpoint of the edge $\bar{V}\bar{W}$ we need only rule out one of
these null elements.

Most of the arguments below depend only on the convex geometry
of the polygons and the Lorentz metric $J$. For such arguments we can apply
arbitrary permutations of the coordinates $x_1, x_2, x_3$, even though
the resulting polygon no longer belongs to the list in Section 3.

\smallskip

\noindent{\bf Trapezoid (Q4):} This polygon is symmetric with respect to
interchanging the first two summands. So we may assume that we are
in case (1) of Table 3 in \cite{DW6} with $V=s=(1, -2, 0), W=u$ and $c=(3\lambda -2, 1-3\lambda, 0)$.
Then $\bar{v}$ corresponds to the (1B) endpoint in \polyc, i.e., $2\bar{v} = \bar{a} + \bar{c}$
with $\bar{a}$ null $\Leftrightarrow  J(\bar{v}, \bar{v}) =J(\bar{c}, \bar{v})$.
Now consider the perspective projection from $\bar{a}$. Since $1<\lambda \leq 3/2$,
and $a=(4-3\lambda, 3\lambda-1, -4)$, the other endpoint of $\Delta^{\bar{a}}$
must correspond to $\bar{u}$. If $J(\bar{u}, \bar{a})=0$, then together with the
null condition for $\bar{c}$ we obtain
$$d_1(d_1+2d_2)^2 =(d_1d_2 -d_1 +4d_2)^2 +d_1d_2(d_1 +6)^2,$$
which has no positive solutions in $d_1$. On the other hand, if $\bar{u}$
corresponds to a (1B) endpoint, then the null vector $2\bar{u}-\bar{a}$
can be connected to $\bar{c}$ by a segment missing \wt. This contradicts Corollary
3.4 in \cite{DW6} and rules out the trapezoid (Q4).

\smallskip

\noindent{\bf Triangle (T11):} First we let $V=s$, $W=u$, so that $s$
is closer to $c$. Since the arguments for (Q4) do not refer to the point $p$,
they apply here. If $V=u$ and $W=s$, the nullity of $2\bar{v}-\bar{c}$
and $1<\lambda$ immediately gives a contradiction.

\smallskip

\noindent{\bf Triangle (T10):} This triangle is symmetric with respect
to interchanging the last two summands. It suffices to consider the
case $V=p, W=w$. Now $\bar{\alpha}$ corresponds to the (1B) endpoint
in \polyc. The null condition for $\bar{a}=2\bar{\alpha}-\bar{c}$
immediately gives a contradiction.

\smallskip

\noindent{\bf Triangle (T9):} First we consider the case $V=w, W=p$.
Then $\bar{v}$ corresponds to the (1B) endpoint in \polyc. The nullity of
the vector $2\bar{v}-\bar{c}$ and $\lambda <1$ immediately gives a contradiction.

For the case $V=p, W=w$, $\bar{v}$ corresponds to
the (1B) endpoint. Besides the nullity of $\bar{c}$ and $\bar{a}=2\bar{v}-\bar{c}$,
we also have  $J(\bar{w}, \bar{a})=0$, for otherwise $\bar{c}$ and $2\bar{w}-\bar{a}$
would be two null vectors in $\mathcal C$ which can be joined by a segment
missing \wt. Putting all these conditions together, we obtain $d_1(d_2 -1) \leq 8d_2,$
so that  $d_1 \leq 16,$ as well as
$$  (d_1 -4)d_2^3 +(d_1^2+ 6d_1 + 4)d_2^2 +d_1(7d_1+25)d_2 + 16d_1^2 =0.$$
One readily verifies that this equation has no positive integer roots.

\smallskip

\noindent{\bf Trapezoid (Q5):} This trapezoid is symmetric with respect to
interchanging the last two summands. So it suffices to consider $V=w, W=p$.
Now $\bar{s}$ corresponds to the (1B) endpoint in \polyc, so
$\bar{a}=2\bar{s}-\bar{c}$ is null. We next consider $\Delta^{\bar{a}}$.
There are now two cases.

Suppose the second endpoint of $\Delta^{\bar{a}}$ corresponds to $\bar{v}$.
If this endpoint is of type (2), then by Tables 1, 2 in \cite{DW6}, we have
$a_1 =4/3,$ contrary to $a_1 =2$. The endpoint cannot be of type (1B)
since our polygon does not contain all three type I vectors.
 Hence the second endpoint is of type (1A), i.e., $J(\bar{a}, \bar{v})=0$,
which implies $d_3 \leq 12$. Together with the nullity of $\bar{c}$ we deduce
$$ (d_3 -4)d_2^3 +(d_3 -2)^2 d_2^2 + 48d_3d_2 + 64d_3^2 = 0.$$
This equation has no positive integer solutions.

The second case occurs if the second endpoint of $\Delta^{\bar{a}}$ corresponds
to $\bar{p}$ and is of type (1A). Up to interchanging the last two coordinates,
we may apply the argument for the second case of (T9).

\smallskip

\noindent{\bf Hexagon (H2) (case (4)):} We consider $c$ collinear with
the edge $vs$. As the hexagon is symmetric with respect to interchanging
the last two coordinates, we only need to examine the case $V=s, W=v$
with $c=(1, -2\lambda, 2\lambda-2).$ By Proposition \ref{c-points}, we
have $1<\lambda<3/2$. ($\lambda=3/2$ would imply that the endpoints of
\polyc are both of type (2), which is not possible.)

Now the type (1B) endpoint in \polyc corresponds to $\bar{w}$, and
$\bar{a}=2\bar{w}-\bar{c}$ is null. Moreover, $a_1 =-1$, so $a$  is
collinear with $x, y$, i.e., the other endpoint of \polya is of type (2).
By Proposition \ref{c-points}, there is a null vector
$\bar{e} \in {\mathcal C}$ lying beyond $\bar{v}$ on the ray from
$\bar{c}$ to $\bar{v}$ and $d_2=d_3$. By symmetry, it follows
that $\bar{p}$ must correspond to the (1B) endpoint in $\Delta^{\bar{e}}$.
So $\bar{b}:=2\bar{p}-\bar{e}=(-1, 4-2\lambda, 2\lambda-4)$ is null.
The nullity of $\bar{a}$ and $d_2=d_3$ give $\lambda=7/6$, which in
turn implies that the nullity of $\bar{e}$ cannot be satisfied.

\smallskip

\noindent{\bf Pentagon (P2) (case (4)):} By Proposition \ref{c-points} it
suffices to consider the case $V=s, W=v$ and $c=(1, -2\lambda, 2\lambda-2)$
with $\lambda \leq 3/2$. The (1B) endpoint in
\polyc corresponds to $\bar{y}$. Let $\bar{a}:=2\bar{y}-\bar{c}$, which
is null. If $\lambda=3/2$, $c=(1, -3, 1),$ and $a=(-3, 1, 1)$ is collinear
with $\bar{x}, \bar{p}$. This configuration is impossible, by Theorem 9.2(i)
of \cite{DW6} (since $r=3$).

So $1< \lambda < 3/2$, which gives $0< a_2 <1$. The other endpoint in \polya
corresponds to $\bar{x}$. If this is of type (1B), we have a null vector
$\bar{b}=2\bar{x}-\bar{a}$. The nullity of $\bar{b}, \bar{c}$ leads to
$$ d_2^2 (d_3-9)+d_3^2(d_2 -9)+d_2d_3(d_2 -13) +d_2d_3(d_3-13)=0.$$
One checks that the only positive integer solutions are
$(d_2, d_3)=(11, 11), (7, 21), (21, 7)$.  Of these only the second is
compatible with $\lambda < 3/2$. But then $d_2$ is not an integer.
So the other endpoint in \polya  must be of type (1A). Now $J(\bar{a}, \bar{x})=0$
and the various null conditions lead to
$$5(4d_2d_3-d_2-d_3)(d_2+d_3)^2 = d_2(2d_2d_3+7d_2-3d_3)^2+d_3(2d_2d_3+17d_2+7d_3)^2,$$
which has no positive integer solutions.

\smallskip

\noindent{\bf Parallelogram (Q3) (case(4)):} We consider $c$ collinear with
the edge $vs$. We have symmetry under interchange of the last two coordinates.
In any case, by Proposition \ref{c-points}, it suffices to consider the case
where $V=s, W=v,$ $c=(1, -2\lambda, 2\lambda-2)$ and $1<\lambda \leq 3/2$.
As in the previous case, we cannot have $\lambda <3/2$.

If $\lambda=3/2$, then $c=(1, -3, 1)$ and $\bar{y}$ corresponds to the type
(1B) endpoint in \polyc. Let $\bar{a}= 2\bar{y}-\bar{c}$, which is null.
The other endpoint in \polya corresponds to $\bar{x}$. It must be of type (1B)
and so $\bar{b}=2\bar{x}-\bar{a}=(1, 1, -3)$ is  null. It follows that
$d_1=d_2=d_3=11$ and (Q3) is inscribed in the triangle $\bar{a}\bar{b}\bar{c}$.
To rule out this special configuration, we apply arguments of the type in the
proof of Theorem 6.4 and in Remark 6.6 of \cite{DW6}. We consider the vertex
$\bar{a}$ and the wedge bounded by the rays $\bar{a}\bar{c}$ and $\bar{a}\bar{b}$.
These arguments show that ${\mathcal C}=\{\bar{a}, \bar{b}, \bar{c}, \bar{z} \}$.
But then the interior points $\beta, \gamma$ correspond to terms unaccounted
for in the scalar curvature function.

\smallskip

\noindent{\bf Hexagon (H2) (case(5)):} Because of the symmetry of the hexagon,
it suffices to consider the case $V=y, W=x$ with $c=(-1, 1-2\lambda, 2\lambda-1)$
and $\lambda \leq 3/2$ (cf Proposition \ref{c-points}).

Note first that $\lambda\neq 3/2$, since we cannot have two type (2) endpoints
in \polyc. So $\lambda < 3/2$, and the (1B) endpoint in \polyc corresponds to $\bar{w}$.
Notice that the null vector $\bar{a}=2\bar{w}-\bar{c}$ has $a_1=1$ and so is collinear
with $\bar{s}, \bar{v}$. Hence the other endpoint in \polya is of type (2), i.e.,
$\bar{a}$ is as in (H2) (case(4)), which has been eliminated.

\smallskip

\noindent{\bf Hexagon (H3):} By Proposition \ref{c-points} we may let $V=y, W=x$
with $c=(-1,1-2\lambda, 2\lambda-1), \lambda < 3/2$. The (1B) endpoint in \polyc
corresponds again to $\bar{w}$. Now the null vector $\bar{a}=2\bar{w}-\bar{c}$ is
collinear with $\bar{z}\bar{v}$. This violates Theorem 9.2(i) of \cite{DW6} as $r=3$.

\smallskip

\noindent{\bf Pentagon (P2) (case (5)):}
By Proposition \ref{c-points} we may assume that $V=y, W=x$ with
$c=(-1, 1-2\lambda, 2\lambda-1), \lambda \leq 3/2$. Then
$\bar{s}$ corresponds to the (1B) endpoint in \polyc. Let $\bar{a}=2\bar{s}-\bar{c}$.
The nullity of $\bar{a}$ gives $d_2=d_1(2\lambda-3)$. The other endpoint of
\polya corresponds to $\bar{v}$. Now $\bar{a}$ cannot be collinear with
$\bar{q}, \bar{v}$ by Theorem 9.2(i) in \cite{DW6}. If $\bar{v}$ corresponds
to a (1B) endpoint in \polya, then $\bar{b}=2\bar{v}-\bar{a}$ is collinear
with $\bar{c}, \bar{x}, \bar{y}$. The nullity of $\bar{b}$ implies that
$\lambda=3/2$, which in turn gives $d_2 =0$.

So $\bar{v}$ corresponds to a (1A) endpoint in \polya. Now $J(\bar{v}, \bar{a})=0$
and the null condition of $\bar{c}$ gives
$$0=(d_3-4)d_2^3 +(d_3+4)(d_3-2)d_2^2 +6d_3(d_3+4)d_2 +36d_3^2.$$
Since $d_3 >1$ we may assume $d_3=2, 3$. The corresponding equation has no
integer roots.

\smallskip

\noindent{\bf Parallelogram (Q3) (case (5)):} The symmetry under interchange
of the last two coordinates means that it suffices to assume that $V=y, W=x$
with $c=(-1, 1-2\lambda, 2\lambda-1)$ and $\lambda \leq 3/2$. The arguments
for (P2) apply here.

\smallskip

\noindent{\bf Pentagon (P5):} This is symmetric under interchange of the
last two coordinates. So we need only consider the case $V=x, W=y$ with
$c=(-1, 2\lambda-1, 1-2\lambda)$ and $\lambda \leq 3/2$. The (1B)
endpoint in \polyc corresponds to $\bar{p}$. Let $\bar{a}=2\bar{p}-\bar{c},$
which is null. If $J(\bar{z}, \bar{a})=0,$ then together with the nullity
of $\bar{c}$ we obtain $8d_2^3 +14d_2^2d_3 +8d_2d_3^2 +2d_3^3 =0$, which
has no positive integer solutions.

Hence $\bar{z}$ corresponds to a (1B) endpoint in \polya. Let $\bar{b}=
2\bar{z}-\bar{a}$. The nullity of $\bar{b}$ gives $d_2=d_3$. So $\lambda=4/3$,
and apart from the fact that $\bar{s}, \bar{v}$ are absent, the configuration
of null vectors in $\mathcal C$ is the same as that in the last
paragraph of (H2)(case (4)). So we can rule it out with the same argument.

\smallskip

\noindent{\bf Pentagon (P3):} The symmetry under interchange of the first
two coordinates implies that we may assume that $c$ is collinear with
$x, y$. Proposition \ref{c-points} implies that we may also assume $V=x, W=y$
with $c=(-1, 2\lambda-1, 1-2\lambda), \lambda<3/2$. The (1B) endpoint of
\polyc corresponds to $\bar{p}$. For the null vector
$\bar{a}=2\bar{p}-\bar{c}$, the other endpoint in \polya is now of type (2).
This contradicts Theorem 9.2(i) in \cite{DW6}.

\smallskip

\noindent{\bf Pentagon (P4):} We may again assume $V=x, W=y$ with
$c=(-1, 2\lambda-1, 1-2\lambda)$ and $\lambda \leq 3/2$. The (1B)
endpoint in \polyc corresponds to $\bar{p}$. But the null vector
$\bar{a}=2\bar{p}-\bar{c}$ is collinear with $\bar{z}, \bar{s}$.
By Theorem 9.2(i) of \cite{DW6} this configuration cannot occur.

\smallskip

\noindent{\bf Trapezoid (Q8):} We may use the arguments for (P5)
above except that now $\bar{b}$ and the null vector beyond $\bar{y}$
(this must be present by Proposition \ref{c-points}) can be joined by
a segment missing \wt, a contradiction.

\smallskip

\noindent{\bf Triangle (T8):} By symmetry we may assume that $V=x, W=y$
and $c=(-1, 2\lambda-1, 1-2\lambda)$ and $\lambda \leq 3/2$. The
(1B) endpoint in \polyc corresponds to $\bar{z}$. Let $\bar{a}=2\bar{z}-\bar{c}$.
Then $\bar{y}$ must correspond to the type (1A) endpoint in \polya. (The
other endpoint in \polya cannot be of type (1B), otherwise we get a segment
joining $\bar{c}$ with the corresponding null vector which would miss \wt.)
The nullity of $\bar{c}$ and $J(\bar{y}, \bar{a})=0$ then lead to
$$0=2d_2d_3^3(2d_2-29)+4d_2^2d_3^2(d_2-2)+2d_2^2d_3(d_2-1) + 54d_3^3 +4d_2d_3^2.$$
Hence $2 \leq d_2 \leq 14,$  and one further finds that $(d_2, d_3)=(3, 1), (7, 7)$
are the only positive integral solutions. In the first case $\lambda =1$, a
contradiction. In the second case, $\lambda=3/2$, but solving for $d_1$
we get $d_1= -d_3$, a contradiction.

\section{\bf Both Endpoints of Type (1B)}

In this section we describe the arguments which will rule out null
elements of $\mathcal C$ for which the endpoints of \polyc are both of
type (1B). We will frequently use arguments from Section 6 of \cite{DW6}.

Our task is quite onerous, as we have to consider the polygons in
Section 3 one by one. In each case, if we replace the edges by the
lines containing them, then the plane is divided into various regions.
Since we have ruled out type (2) endpoints, our null vector $\bar{c}$
must lie in the interior of one of the regions. The assumption on $\bar{c}$
implies that there are vertices $\bar{V}, \bar{W}$ of the polygon such
that $\bar{a}:= 2\bar{V}-\bar{c}$ and  $ \ap := 2\bar{W}-\bar{c}$
are null elements of $\mathcal C$. We will use this notation throughout this
section.

A number of the regions can be ruled out immediately by the fact that
$\bar{a} \ap$ misses \wt. These are typically regions whose closure
intersects a single vertex of \wt. Recall as well that the methods in
\cite{DW6} already allow us to rule out those polygons which do not contain
all three type I vectors (cf Theorem 8.7 and remarks before Theorem 6.12
in \cite{DW6}).

\smallskip

Let us consider first the triangles in Section 3. Triangle (T1) is course
possible, being the situation of a triple warped product, i.e., case (1)
of Theorem \ref{smallrthm}. The next four triangles do not contain all three
type I vectors and so have been eliminated. (T6)-(T8) are triangles containing
at least one member of \wt as the midpoint of an edge. These are eliminated by
Remark 6.13 of \cite{DW6}. (T9)-(T10) have exactly one edge with interior points
in \wt and no interior points in \wt. They are ruled out by Remark 6.14 of \cite{DW6}.

We are left with (T11) and (T12), both of which contain an interior point
in \wt. We discuss the first case in detail and leave the second case to the
reader.

\smallskip

\noindent{\bf Triangle (T11):} First we suppose that $\bar{c}$ lies in the
region of the plane sharing a boundary with $\bar{u}\bar{v}$. Let
$\bar{a}=2\bar{u}-\bar{c}$ and $\ap=2\bar{v}-\bar{c}$. By arguments similar
to those on pp. 616-617 of \cite{DW6}, we conclude that all elements of
${\mathcal C} \setminus \{\bar{c}, \bar{a}, \ap \}$ must lie in the half-plane
bounded by $\bar{a}\ap$ on the opposite side of $\bar{c}$. If
$\bar{s} \notin {\mathcal C},$ then $\bar{a}, \bar{s}$ and $\ap$ must be
collinear and ${\mathcal C} = \{\bar{c}, \bar{a}, \ap \}.$ But the term in
the superpotential equation corresponding to $\bar{\alpha}$ would be
unaccounted for. Thus $\bar{s} \in {\mathcal C}$ and we have $J(\bar{a}, \bar{s})=0
=J(\bar{s}, \ap).$ When we write these out in detail, they turn out to be
inconsistent with each other.

The case where $\bar{c}$ lies in the region of the plane sharing a boundary
with $\bar{u}\bar{s}$ is handled in an analogous manner.

Finally, consider the case where $\bar{c}$ lies in the region of the plane
sharing a boundary with $\bar{u}\bar{s}.$ In this case, all elements of
${\mathcal C} \setminus \{\bar{c}, \bar{a}, \ap, \frac{1}{3}(2\ap + \bar{a}),
\frac{1}{3}(\ap +2\bar{a}) \}$ must lie in the half-plane bounded by
$\bar{a}\ap$ on the opposite side of $\bar{c}.$ We must now have $\bar{v} \in
{\mathcal C}$ and $J(\bar{a}, \bar{v})=0=J(\bar{v}, \ap).$ These again turn out
to be inconsistent with each other.

\smallskip

Among the quadrilaterals in Section 3, we consider the following examples
which provide a sample of the arguments used.

\noindent{\bf Square with midpoint (Q2):} We note that the scalar curvature
function is realised by $K=K_1 \times K_2 \subset H_1 \times H_2 \subset G,$
where $H_i/K_i$ and $G/(H_1 \times H_2)$ are isotropy irreducible. The following
diagram will be useful in the discussion below.

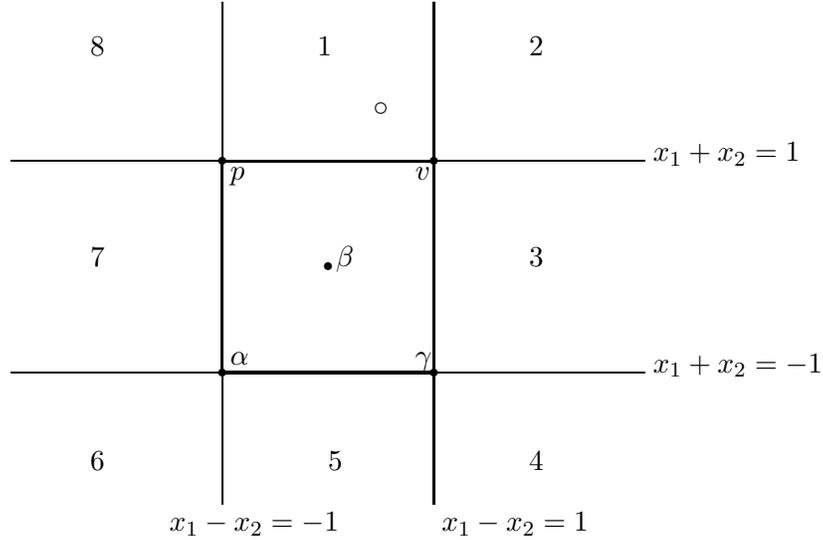
\begin{picture}(200, 210)(10, 0)
\put(120,140){\line(1,0){240}}
\put(200, 140){\circle*{3}}
\put(280, 140){\circle*{3}}
\put(200, 60){\circle*{3}}
\put(280, 60){\circle*{3}}
\put(120, 60){\line(1, 0){240}}
\put(200, 10){\line(0, 1){190}}
\put(280, 10){\line(0, 1){190}}
\put(240, 100){\circle*{3}}
\put(363, 140){$x_1+x_2 =1$}
\put(363, 60){$x_1+x_2 =-1$}
\put(283, 0){$x_1 -x_2 =1$}
\put(180, 0){$x_1-x_2 =-1$}
\put(243, 100){$\beta$}
\put(203, 133){$p$}
\put(273, 133){$v$}
\put(203, 63){$\alpha$}
\put(273, 63){$\gamma$}
\put(236, 180){$1$}
\put(316, 180){$2$}
\put(316, 100){$3$}
\put(316, 23){$4$}
\put(240, 23){$5$}
\put(150, 23){$6$}
\put(150, 100){$7$}
\put(150, 180){$8$}
\put(260, 160){\circle{4}}
\thicklines
\put(200, 140){\line(1,0){80}}
\put(200, 60){\line(1,0){80}}
\put(200, 60){\line(0,1){80}}
\put(280, 60){\line(0,1){80}}
\end{picture}

\smallskip

Note that interchanging the first two coordinates is a symmetry for the
square. So we need only consider $\bar{c}$ lying in regions $1$ to $5$.
To ensure that the null vectors $\bar{a}$ and $\ap$ cannot be connected
by a segment which avoids \wt, $\bar{c}$ must further be in regions $1, 3$
or $5$.

\noindent{\bf A. $\bar{c}$ in region $1$:}  We have $-1<c_1 -c_2 <1,\, c_1 +c_2 >1$.

The vector $a=2p-c=(-c_1, 2-c_2,-4-c_3)$ lies in region $7$ since
$a_1+a_2 =2-c_1-c_2 <1$ and if $a_1+a_2 <-1$, then the segment $\bar{a}\ap$
would miss \wt. Likewise, $a^{\prime}=2v-c=(2-c_1, -c_2, -4-c_3)$ must
lie in region $3$. (Since we have ruled out the presence of type (2)
endpoints we may assume that $\bar{a}$ and $\ap$ are not collinear with
$\bar{\alpha}, \bar{\gamma}$.) As observed in \cite{DW6}, the null condition of $\bar{a}$, given
the nullity of $\bar{c},$ is equivalent to $J(\bar{p}, \bar{p})=J(\bar{p}, \bar{c}),$
which is $\frac{1}{d_2}+\frac{4}{d_3}=\frac{c_2}{d_2}-\frac{2c_3}{d_3}.$
Similarly, the nullity of $\ap$ is equivalent to $J(\bar{v}, \bar{v})=J(\bar{v}, \bar{c}),$
which is $\frac{1}{d_1}+\frac{4}{d_3}=\frac{c_1}{d_1}-\frac{2c_3}{d_3}.$

Next we consider the wedge bounded by the rays $\bar{c}\bar{a}$ and $\bar{c}\ap$.
Using arguments similar to those in the proof of Theorem 6.4 in \cite{DW6}, we
see that all elements of ${\mathcal C}\setminus \{\bar{c}, \bar{a}, \ap \}$
lie below the line $\bar{a}\ap$. Let $\bar{e}$ be the intersection of
$\bar{a}\bar{\alpha}$ and $\ap\bar{\gamma}$. There are now two possibilities.

If $\bar{e} \notin {\mathcal C},$ then we must have $\bar{\alpha}, \bar{\gamma}
\in {\mathcal C}$ and $J(\bar{a}, \bar{\alpha})=0=J(\bar{\gamma}, \ap).$
But by Proposition 3.7 in \cite{DW6}, we have $J(\bar{\alpha}, \bar{\gamma})=0$
as well. But this is not satisfied.

Hence $\bar{e} \in {\mathcal C}$ and is null. Furthermore, we must have
$\bar{e}=2\bar{\alpha}-\bar{a} =2\bar{\gamma}-\ap$ and ${\mathcal C}=
\{\bar{a}, \bar{c}, \bar{e}, \ap \}.$  The null condition $J(\bar{\alpha}, \bar{\alpha})
=J(\bar{\alpha}, \bar{a})$ for $\bar{e}$ leads to $c_1=1$. Using the
other expression for $\bar{e}$ we obtain $c_2=1,$ which is expected from
symmetry. So $c_3 =-3$ since $c_1+c_2+c_3 =-1.$ But this contradicts the
null condition for $\bar{a}$.

\noindent{\bf B. $\bar{c}$ in region $3$:}  $c_1-c_2 > 1, \, -1< c_1+c_2 <1$

The arguments are similar to those above. All elements of ${\mathcal C} \setminus
\{\bar{c}, \bar{a}, \ap \}$ must lie to the left of the line $\bar{a}\ap.$
Let $\bar{e}$ be the intersection of $\bar{a}\bar{p}$ and $\ap\bar{\alpha}$.
Then $J(\bar{p},\bar{\alpha})=\frac{1}{4} \neq 0$ means that $\bar{e}$ is null
and equal to both $2\bar{p}-\bar{a}$ and $2\bar{\alpha}-\ap.$ The null condition
for $\ap$ gives $c_2=-1$ while the null condition $J(\bar{\alpha}, \bar{\alpha})=
J(\bar{\alpha}, \ap)$ for $\bar{e}$ gives $c_1=1$. Hence $c_3 =-1$. Putting these
in the null condition $J(\bar{p}, \bar{p})=J(\bar{p}, \bar{a})$ (also for $\bar{e}$)
gives a contradiction.

\noindent{\bf C. $\bar{c}$ in region $5$:} $-1<c_1-c_2 <1, \, c_1+c_2 <-1$

We may assume that we are not in the configuration of subcase A above.
In other words, the intersection $\bar{e}$ of the lines $\ap \bar{v}$ and
$\bar{a}\bar{p}$ is not null. Furthermore, $\bar{p}, \bar{c} \in {\mathcal C}$
and $J(\bar{a}, \bar{p})=0=J(\ap, \bar{v}).$ This time, applying Proposition
3.7 in \cite{DW6} gives $d_3 =4$ instead of an immediate contradiction.
(This is indicative of the fact that while the square has reflection
symmetry about the perpendicular bisector of $pv$, it is not induced by
permutations of the summands.) However, the null conditions for $\bar{a}$
and $\ap$ give $c_1=c_2=-1, c_3=1$. The orthogonality conditions
above imply that $d_1=d_2=2$. These values now contradict the null condition
for $\bar{c}$.

\hspace{0.2cm}

In analysing adjacent (1B) vertices in \polyc in \cite{DW6} we did not have to
deal with pentagons. We next discuss the case of (P3) in detail as an example
of how (P1)-(P5) are eliminated.

\noindent{\bf Pentagon (P3):} The pentagon is symmetric under the interchange
of the first two coordinates. We will refer to the following (schematic) diagram.

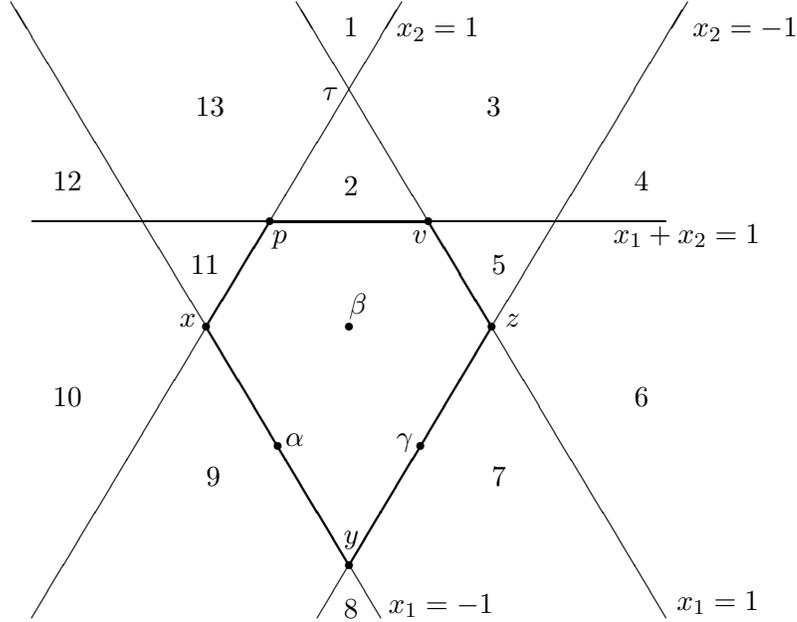
\begin{picture}(200, 240)(10, 0)
\put(360, 140){\line(-1,0){240}}
\put(270, 140){\circle*{3}}
\put(210, 140){\circle*{3}}
\put(186, 100){\circle*{3}}
\put(294, 100){\circle*{3}}
\put(240, 100){\circle*{3}}
\put(120, -10){\line(3,5){140}}
\put(240, 10){\circle*{3}}
\put(228, -10){\line(3,5){140}}
\put(252, -10){\line(-3,5){140}}
\put(360, -10){\line(-3,5){140}}
\put(213, 55){\circle*{3}}
\put(267, 55){\circle*{3}}
\put(238, 210){$1$}
\put(238, 150){$2$}
\put(182, 180){$13$}
\put(292, 180){$3$}
\put(211, 132){$p$}
\put(264, 132){$v$}
\put(128, 152){$12$}
\put(348, 152){$4$}
\put(176, 100){$x$}
\put(299, 100){$z$}
\put(180, 120){$11$}
\put(294, 120){$5$}
\put(240, 105){$\beta$}
\put(238, 18){$y$}
\put(128, 70){$10$}
\put(348, 70){$6$}
\put(186, 40){$9$}
\put(294, 40){$7$}
\put(238, -10){$8$}
\put(255, -8){$x_1 =-1$}
\put(364, -7){$x_1=1$}
\put(340, 132){$x_1+x_2 =1$}
\put(370, 210){$x_2 =-1$}
\put(258, 210){$x_2=1$}
\put(230, 186){$\tau$}
\put(216, 55){$\alpha$}
\put(258, 55){$\gamma$}
\thicklines
\put(186, 100){\line(3,5){24}}
\put(210, 140){\line(1,0){60}}
\put(294, 100){\line(-3, 5){24}}
\put(240, 10){\line(3,5){54}}
\put(240, 10){\line(-3,5){54}}
\end{picture}

\hspace{1.5cm}

A null $\bar{c} \in {\mathcal C}$ such that both endpoints of \polyc
are of type (1B) must lie in one of the regions labelled 1 through 13.
By symmetry, we need only consider those in regions 1 through 8. We can
eliminate regions 1, 4, 8 since the segment connecting the null vectors
$\bar{a}$ and $\ap$ actually misses \wt. Note also that the triangle
$\bar{x}\bar{y}\bar{z}$ is equilateral, the angles at $\bar{x}$ and $\bar{z}$
are both $2\pi/3$, and the lines $x_2=1$ and $x_1=1$ intersect at
$\tau=(1, 1, -3)$.

\noindent{\bf A. $\bar{c}$ in region 3:} $-1<c_2<1, c_1+c_2 > 1, c_1 >1$

The null conditions for $\bar{a}=2\bar{p}-\bar{c}$ and $\ap=2\bar{z}-\bar{c}$
are respectively $\frac{1}{d_2}+\frac{4}{d_3}= \frac{c_2}{d_2}-\frac{2c_3}{d_3}$
and $\frac{1}{d_1}+\frac{1}{d_2} +\frac{1}{d_3} = \frac{c_1}{d_1}-\frac{c_2}{d_2}
-\frac{c_3}{d_3}.$ One checks that $\bar{a}$ must lie in region $10$, while $\ap$
lies in region $8$ or $7$. (It cannot lie on the line $\bar{x}\bar{y}$ as we have
ruled out type (2) endpoints.) In the former case, the segment $\bar{a}\ap$
would miss \wt. So $\ap$ is in region $7$, which means that $\ap_1 =2-c_1 > -1.$

Next we consider the null element $\ap$. If $\bar{y} \notin {\mathcal C}$, then
the only possibility is for $2\bar{y} = \bar{a} +\ap.$ Hence $c=(2, 1, -4)$
and this contradicts the null condition for $\bar{a}.$ If $\bar{y} \in {\mathcal C}$
then we have in addition $J(\bar{a}, \bar{y})=0=J(\ap, \bar{y})$. From this we
obtain $\frac{1}{d_1}=\frac{2}{d_2} +\frac{1}{d_3},$ as well as
$$ \frac{c_1}{d_1}=\frac{3}{4} +\frac{7}{4d_2} + \frac{3}{2d_3}, \,\,\,
   \frac{c_2}{d_2}= \frac{1}{2}\left(1-\frac{1}{d_2}+\frac{2}{d_3} \right), \,\, \,
   \frac{c_3}{d_3}=\frac{1}{2}\left(\frac{1}{2}-\frac{3}{2d_2} -\frac{3}{d_3} \right). $$
Substituting these expressions into the null condition for $\bar{c}$ and
simplifying, one obtains $0=d_2d_3^4(d_2 -6) +9d_2^2d_3^3(d_2-2) + \mbox{{\rm positive
 quantity}}$. So $1 \leq d_2 \leq 5.$ For each of these values, the resulting polynomial
in $d_3$ has no integral roots.

\noindent{\bf B. $\bar{c}$ in region $6$:} $c_1+c_2 <1,\, c_2 < -1, \,c_1 >1$

The null condition for $\bar{a}=2\bar{v}-\bar{c}$ is $\frac{1}{d_1}+\frac{4}{d_3}=
\frac{c_1}{d_1} -\frac{2c_3}{d_3}$ and that for $\ap=2\bar{y}-\bar{c}$ is
$\frac{1}{d_1}+\frac{1}{d_2} +\frac{1}{d_3} = -\frac{c_1}{d_1}-\frac{c_2}{d_2} +
\frac{c_3}{d_3}.$ It follows that $\ap$ lies in region $9$. Note also that
as $a_1^{\prime}< -3, a_1^{\prime}+a_2^{\prime} > -5,$ we have $a_2 =-c_2 < 3,$
and we may assume that $a_1 > 1$ in order that $\bar{a}\ap$ would not miss \wt.

We next consider $\bar{x}$. If $\bar{x} \notin {\mathcal C}$ then $2\bar{x}=
\bar{a}+\ap$, so that $c=(1, -2, 0)$. This contradicts the nullity of $\bar{a}.$
So $\bar{x} \in {\mathcal C}$ and $J(\bar{a}, \bar{x})=0=J(\bar{x}, \ap).$
Subtracting these conditions yields $\frac{2}{d_1}=\frac{1}{d_2} + \frac{3}{d_3}.$
Solving for $\frac{c_i}{d_i}$ we get
$$ \frac{c_1}{d_1}= \frac{1}{2}\left(1-\frac{1}{d_1}+\frac{1}{d_2}-\frac{5}{d_3} \right), \,\,\,
   \frac{c_2}{d_2}=-\frac{1}{4}\left(1+\frac{13}{d_1}-\frac{3}{d_2} + \frac{7}{d_3} \right), \,\,\,
   \frac{c_3}{d_3}=\frac{1}{4}\left(1-\frac{3}{d_1}+\frac{1}{d_2} -\frac{13}{d_3} \right).$$
Substituting these into the null condition for $\bar{c}$ and using the relation
between the dimensions, one obtains in the case when $d_3 \geq 15$ that
\begin{eqnarray*}
0&=&(4d_2^2-4d_2+1)d_3^4 + (48d_2^2 -156d_2+130)d_2d_3^3 +
    (12d_2^2-732d_2+2116)d_2^2d_3^2 \\
  & & +(636d_2+10062)d_2^3d_3 + 8427d_2^4  \\
  &>& d_2d_3^2\left( 15(48d_2^2-156d_2 +130) +d_2(12d_2^2-732d_2 +2116)\right) \\
  &=& d_2d_3^2(12d_2^3 -12d_2^2 -224d_2 + 1950),
\end{eqnarray*}
where we have also used the fact that the coefficient of the $d_3^4$ term is always
positive. For $d_2>0$ one easily verifies (by finding the minimum) that the last expression
is positive. Hence we may assume that $2 \leq d_3 \leq 14.$ For these values
of $d_3$ one checks that the right-hand side of the first equation above is
always positive.

\noindent{\bf C. $\bar{c}$ in region $7$:}  $-1<c_1<1, \,\, c_2<-1$

By subtracting the null conditions for $\bar{a}=2\bar{z}-\bar{c}$ and
$\ap=2\bar{y}-\bar{c},$ we obtain $\frac{c_2}{d_2}=-\frac{1}{d_1}-\frac{1}{d_2}
-\frac{1}{d_3}$ and $\frac{c_1}{d_1}=\frac{c_3}{d_3}.$ In particular $c_1=\frac{d_2}{d_3}.$
Note that $\ap$ must lie in region $9$ since we have already eliminated type (2)
endpoints and $\bar{a}\ap$ must meet \wt. The null condition for $\bar{c}$
is given by
$$0=-d_1^2d_2d_3^2 +d_1^2(d_2+d_3)^2 +d_1d_2^3 +2d_1d_2d_3(d_2+d_3) +d_2^2d_3(d_2+d_3).$$
There are now two possibilities.

If $\bar{x} \in {\mathcal C}$ then $J(\ap, \bar{x})=0.$ It follows from this that
$d_1(d_2 +d_3 +d_2d_3) =d_2(2d_2 +3d_3).$ Substituting the value of $d_1$ from
this equation into the null condition for $\bar{c}$ and using $d_2^2 \geq d_2,$ we
obtain a positive expression, a contradiction.

If $\bar{x} \notin {\mathcal C}$, then the vector $\bar{e}^{\prime} = 2\bar{x} -\ap$
is null. This null condition implies that $d_1 =d_2, c_3 =1,$ so $c_1+c_2 =-2$ and
$e_1^{\prime}+ e_2^{\prime} = 2 = a_1 +a_2.$ This means that $\bar{a}\bar{e}^{\prime}$
misses \wt, a contradiction.

\noindent{\bf D. $\bar{c}$ in region $5$:} $c_2>-1, \, c_1>1, \, c_1+c_2 < 1$

The null condition for $\bar{a}=2\bar{v}-\bar{c}$ is $\frac{1}{d_1}+\frac{4}{d_3}=
\frac{c_1}{d_1} - \frac{2c_3}{d_3},$ and that for $\ap=2\bar{z}-\bar{c}$ is
$\frac{1}{d_1}+\frac{1}{d_2}+\frac{1}{d_3} =\frac{c_1}{d_1}-\frac{c_2}{d_2}-\frac{c_3}{d_3}.$
One checks using the inequalities for $c$ that $\ap$ lies in region $7$. By the
previous case, we must have $\bar{y} \in {\mathcal C}$ and $J(\ap, \bar{y})=0.$
We can then solve for $\frac{c_i}{d_i}$ to get
$$ \frac{c_1}{d_1}=1+\frac{2}{d_1}-\frac{1}{d_2}-\frac{1}{d_3},\,\,\,
   \frac{c_2}{d_2}=\frac{1}{2} \left( 1+\frac{1}{d_1}-\frac{3}{d_2}+\frac{1}{d_3} \right), \,\,\,
  \frac{c_3}{d_3}=\frac{1}{2}\left(1+\frac{1}{d_1}-\frac{1}{d_2}-\frac{5}{d_3} \right). $$
The null condition for $\bar{c}$ becomes
$$0=d_1^2d_2(d_2-2)d_3^3 +d_1^2d_2^2(d_2-4)d_3^2 + 2d_1d_2(d_2-1)d_3^3 +
    2d_1^3d_2d_3(d_2(d_3-4)+d_3(d_2-4))+ \mbox{\rm pos. quantity}.$$
Hence either $d_2 \leq 3$ or $2 \leq d_3 \leq 3$. Substituting in these values,
we can then verify that the resulting polynomials do not admit positive
integral solutions.

\noindent{\bf E. $\bar{c}$ in region $2$:}  $c_2 <1,\, c_1<1, \, c_1+c_2 > 1$

The null condition for $\bar{a}=2\bar{p}-\bar{c}$ is $\frac{1}{d_2}+\frac{4}{d_3}=
\frac{c_2}{d_2}-\frac{2c_3}{c_3}$ and that for $\ap=2\bar{v}-\bar{c}$ is
$\frac{1}{d_1}+\frac{4}{d_3}= \frac{c_1}{d_1}-\frac{2c_3}{c_3}$. By the previous
case we must have $\bar{z} \in {\mathcal C}$ and $J(\ap, \bar{z})=0.$ We now
get
$$ \frac{c_1}{d_1}=2-\frac{1}{d_1}-\frac{2}{d_2}-\frac{4}{d_3}, \,\,\,
   \frac{c_2}{d_2}= 2-\frac{2}{d_1}-\frac{1}{d_2}-\frac{4}{d_3},\,\,\,
   \frac{c_3}{d_3}= 1-\frac{1}{d_1}-\frac{1}{d_2}-\frac{4}{d_3}. $$
Using these we can write the condition $c_1+c_2+c_3 =-1$ in the form
\begin{eqnarray*}
0 &=& \frac{1}{2}d_1d_2d_3(d_1-5)+ \frac{1}{2}d_1d_2d_3(d_2 -5) +d_1d_3^2(\frac{d_2}{2}-1)
       +d_2d_3^2(\frac{d_1}{2}-1) + \\
  & & d_1^2d_3(\frac{d_2}{2}-2)+d_2^2d_3(\frac{d_1}{2}-2) +d_1d_2^2(d_3-4)+d_1^2d_2(d_3-4).
\end{eqnarray*}
Hence one of $d_1 \leq 4, d_2 \leq 4$ or $2 \leq d_3 \leq 3$ must be true. We now
use these values in the null condition for $\bar{c}$. By symmetry, the cases
$d_1 \leq 4$ and $d_2 \leq 4$ are analogous. If we let $d_1=1$ or$d_3 =2$ we
immediately obtain positive expressions for the null condition. If $d_3 =3,$
the null condition gives an equation for $d_1, d_2$ which requires $d_1 \leq 11$
(or $d_2 \leq 11$ by symmetry). For these values of $d_1$, the resulting polynomial
in $d_2$ has no integer roots. Similar arguments rule out the  $d_1=2,3,4$
cases.

\smallskip

The cases of the hexagons (H2) and (H3) are quite analogous, except that there
are more regions for $\bar{c}$ to lie in. In the case of (H3), which is not symmetric,
there are $19$ regions to consider. Apart from the tedium, there are no new ideas
necessary for ruling out all the cases.

\section{\bf Case of Endpoints of Type (1A)(1B)}

We have thus far shown that for any null $\bar{c} \in {\mathcal C}$ (which does
not correspond to a type I vector), one of the endpoints in \polyc is of
type (1A) and the other must be of type (1B), unless we are in the situation of
case (1) of Theorem \ref{smallrthm}. In this section we outline the
arguments which show that the (1A)(1B) combination cannot occur.

We begin by describing in detail the arguments for eliminating parallelogram (Q1).
Recall that by Theorem 3.14 in \cite{DW6}, we may assume that none of the null
elements of $\mathcal C$ is of type I.

\smallskip

\noindent{\bf Parallelogram (Q1):} Let $\bar{c}, \bar{a}$ denote null elements
of $\mathcal C$. There are four cases, corresponding to $\bar{c}+\bar{a}=
2\bar{v}, 2\bar{s}, 2\bar{p}$ and $2\bar{\gamma}$. For each case, we consider
the conditions imposed by the (1A) endpoint in \polyc and in \polya. Together
with the nullity of $\bar{c}$ and $\bar{a}$, one derives certain equations,
some of which are diophantine in nature, which have no admissible solutions.

\noindent{\bf Subcase I: $\bar{c}+\bar{a} = 2\bar{v}$}

Let $c=(c_1, c_2, c_3)$, then $a=(2-c_1, -c_2, -4-c_3)$. Given that $\bar{c}$ is null,
the null condition for $\bar{a}$ is $J(\bar{v}, \bar{v})=J(\bar{c}, \bar{v})$,
which is $\frac{1}{d_1} + \frac{4}{d_3} = \frac{c_1}{d_1}-\frac{2c_3}{d_3}$.

If $\bar{\gamma}$ is (1A) for both $\bar{c}$ and $\bar{a}$ (this happens when
$\bar{c}\bar{a}$ is very long), then $J(\bar{v}, \bar{\gamma})=0$ (since
$\bar{v}$ is the midpoint of $\bar{a}\bar{c}$), which does not hold. This
observation can be used in many of the cases we have to consider.

Another general observation can be used to show that $\bar{p}$
cannot be (1A) for both $\bar{c}$ and $\bar{a}$. Consider the
pair of rays $\bar{p}\bar{\gamma}$ and $\bar{p}\bar{v}$. If $\bar{p}$
corresponds to a (1A) endpoint of $\Delta^{\bar{b}}$ for some null element
$\bar{b} \in \mathcal C$, then $\bar{b}$ must lie in one of the two open
``quadrants" bounded by exactly one of the two rays. Since $vs$ is parallel
to $p\gamma$, $\bar{c}$ and $\bar{a}$ cannot both satisfy the above condition.
Likewise, $\bar{s}$ cannot be (1A) for $\bar{c}$ and $\bar{a}$.

Observe next that we cannot have $J(\bar{c}, \bar{\gamma})=0$ (resp.
$J(\bar{a}, \bar{\gamma})=0$) since this would force  $c$ (resp. $a$) also to
be a type I vector. So we end up with $J(\bar{c}, \bar{p})=0$
and $J(\bar{a}, \bar{s})=0$, with $\bar{p}, \bar{s} \in {\mathcal C}$.
These and the null condition for $\bar{a}$ give
$$\frac{c_1}{d_1}= \frac{8}{d_3}-1,\,\, \frac{c_2}{d_2}= \frac{4}{d_3}-\frac{1}{d_1},
  \,\,\frac{c_3}{d_3}=\frac{1}{2}\left(-1-\frac{1}{d_1} + \frac{4}{d_3}\right).$$
Using the relation $c_1+ c_2 + c_3 =-1$ we see that
$$ 2d_2(4d_1 -d_3)= -16d_1^2 +2d_1^2 d_3 -6d_1d_3 +d_1 d_3^2 + d_3^2.$$
Now from the null condition for $\bar{c}$ it follows that $4d_1 \neq d_3$.
So the above equation can be used to express $d_2$ in terms of $d_1, d_3$
in the null condition for $\bar{c}$. We then obtain
$$0=d_3^3(d_1^2 -1)+ d_1^2d_3^2(d_1 -8) + 3d_1^3d_3(d_3-16) + \,\, \mbox{\rm positive quantity}.$$
So either $d_1 \leq 7$ or $d_3 \leq 15$. For these values, one easily verifies
(e.g. using MAPLE) that the only integral solution of the last
equation is $(d_1, d_3)=(3, 6)$. But then $d_2 =0$, a contradiction.

\noindent{\bf Subcase II: $\bar{c}+\bar{a} = 2\bar{s}$}

As in the previous subcase, we may assume that $J(\bar{c}, \bar{\gamma}) \neq 0$
and $J(\bar{a}, \bar{\gamma}) \neq 0$. Also, $\bar{v}$ cannot be (1A) for both
$\bar{c}$ and $\bar{a}$. Since $J(\bar{p}, \bar{s})=1 + \frac{2}{d_2} > 0$,
$\bar{p}$ cannot be (1A) for $\bar{c}$ and $\bar{a}$. Finally, if
$J(\bar{c}, \bar{v}) = 0=J(\bar{a}, \bar{p})$, then it follows that
$$ \frac{c_1}{d_1}=4-\frac{1}{d_1} +\frac{4}{d_2} > 3, \,\, \frac{c_2}{d_2}= 2-\frac{1}{d_1},
 \,\,\, \frac{c_3}{d_3} = \frac{3}{2} -\frac{1}{2d_2} + \frac{2}{d_2}. $$
These relations imply that the null condition for $\bar{c}$ is violated.

\noindent{\bf Subcase III: $\bar{c} + \bar{a} = 2\bar{p}$}

As in the previous two subcases, we need only to consider the situation
$J(\bar{c}, \bar{v})=0=J(\bar{a}, \bar{s})$. After solving for $c_i/d_i$,
the condition $c_1 + c_2 +c_3 = -1$ becomes
$$0=d_1d_2(3d_3-8) + 2d_2^2(d_3 -2) + \, \mbox{\rm positive quantity},$$
which implies that $d_3 =2$. Substituting this into the null condition
for $\bar{c}$, we obtain the contradiction
$$ 0= d_1 (d_2 -2)^2 + \,\, \mbox{\rm positive quantity} > 0$$

\noindent{\bf Subcase IV: $\bar{a}+ \bar{c} = \bar{\gamma}$}

In this case, the nullity of $\bar{c}$ yields $c_2 = -1,$ so $c_1 + c_3 = 0.$
As before we quickly see that none of $\bar{v}, \bar{p},$ or $\bar{s}$ can
be (1A) for both $\bar{c}$ and $\bar{a}$. This leaves us with the three cases
(i) $J(\bar{c}, \bar{p})=0=J(\bar{a}, \bar{s})$, (ii) $J(\bar{c}, \bar{v})=0=J(\bar{a}, \bar{s})$,
and (iii) $J(\bar{c}, \bar{p})=0=J(\bar{a}, \bar{v}).$ In the first two cases, after
solving for $c_1$ and $c_3$, we find that $c_1 +c_3 > 0$. In the third case,
after solving for $c_1$ and $c_3$, we see that the null condition for $\bar{c}$ is
violated.

Therefore, we have ruled out parallelogram (Q1).

\smallskip

All the other polygons in Section 3 are eliminated by similar arguments.
In the following, using the same notation as above, we will highlight those
cases which present additional difficulties.

\smallskip

\noindent{\bf Triangle (T2):} When $\bar{c}+\bar{a} = 2\bar{\beta}$ with
$J(\bar{c}, \bar{v})=0=J(\bar{a}, \bar{p})$, we obtain, using $c_1+c_2+c_3 =-1$
and the assumption that $c$ is not type I, $c_3=-1, d_1=d_2$. The null
condition for $\bar{c}$ becomes $1=2d_1(1-\frac{2}{d_3})^2 + \frac{1}{d_3}$.
The only positive integral solution is  $d_1=d_2=d_3 =3,$ and hence
$c=(1, -1, -1), a=(-1, 1, -1).$ Let $\bar{c}\bar{v}$ and $\bar{a}\bar{p}$ intersect
at $\bar{b}$. We have (T2) inscribed in the triangle $\bar{c}\bar{a}\bar{b}$
such that the vertices of (T2) are the midpoints of its sides.
Note that all elements of $\mathcal C$ must lie in $\bar{c}\bar{a}\bar{b}$. Consider
the wedge bounded by the rays $\bar{c}\bar{b}$ and $\bar{c}\bar{a}$. All
elements of $\mathcal C$ lying in the wedge but not on $\bar{c}^{\perp}$
must have a positive inner product with $\bar{c}$. There are also no elements
of $\mathcal C$ in the triangle $\bar{c}\bar{\beta}\bar{v}$ apart from
$\bar{c}$ and $\bar{v}$. It follows that the term in the superpotential
equation corresponding to $\bar{c}+\bar{p}$ is unaccounted for.

\smallskip

\noindent{\bf Triangle (T9):} If we have $\bar{c}+\bar{a} =2\bar{p}$ with
$J(\bar{c}, \bar{w})=0=J(\bar{a}, \bar{v})$, then after solving for
$c_i$ and using $c_1+c_2 +c_3 =-1$, we obtain (after some simplification)
$$ 0=d_1d_2(5d_3-8)+d_2d_3(2d_2+d_3+6) +4d_2^2  + 2d_3^2 + 4d_1d_3,$$
which cannot hold since $d_3 \geq 2$.

If $\bar{c}+\bar{a} = 2\bar{v}$ with $J(\bar{c}, \bar{p})=0=J(\bar{a}, \bar{w}),$
then the condition $c_1+c_2+c_3 =-1$ yields $0=10d_3 +20d_1 +8d_2 +d_2d_3 -2d_1d_3 -d_3^2$.
This can be rewritten as $d_2(d_3+8)=(2d_1+d_3)(d_3-10)$, so $d_3>10$.
The null condition for $\bar{c}$ gives (after some simplification)
$$2d_1d_3^2(d_3-\frac{31}{2})  + 6d_1^2 d_3(d_3-14) +240 d_1^2 +56d_1d_3 +9d_3^2 =0,$$
where we have used the previous equation to eliminate $d_2$. It follows that
$d_3 < 16$. Substituting the values of $d_3$ between $11$ and $15$ now give
equations in $d_1$ with no integer roots.

\smallskip

\noindent{\bf Triangle (T10):} It is interesting to note that $\mathcal W$
is actually realised by the scalar curvature function of $(G_1/K_1) \times (G_2/K_2)$
where $G_1/K_1$ is isotropy irreducible and $\kf_2$ is a maximal
${\rm Ad}_{K_2}$-invariant subalgebra of $\g_2$.

Here the subcase $\bar{a}+\bar{c}=2\bar{\alpha}$ with
$J(\bar{c}, \bar{p})=0=J(\bar{a}, \bar{w})$ requires special attention.
Since $c, a$ are assumed not to be of type I, one quickly gets $c_1=-1$, $d_2 =d_3$,
and $c_2=-c_3=d_2/3$. Applying the null condition for $\bar{c}$ one immediately
gets $d=(3,3,3)$ or $d=(9, 4, 4)$. In the first case, $c=(-1, 1, -1),  a=(-1, -1, 1)$
while in the second case $c=(-1, \frac{4}{3}, -\frac{4}{3}), a=(-1, -\frac{4}{3}, \frac{4}{3}).$
The following diagram represents the configuration to be analysed.

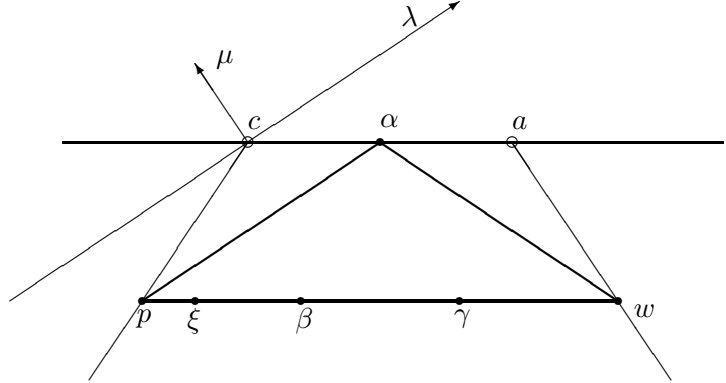
\begin{picture}(200, 150) (10, 10)
\put(150, 40){\circle*{3}}
\put(210, 40){\circle*{3}}
\put(270, 40){\circle*{3}}
\put(330, 40){\circle*{3}}
\put(240, 100){\circle*{3}}
\put(190, 100){\circle{4}}
\put(290, 100){\circle{4}}
\put(190, 105){$c$}
\put(290, 105){$a$}
\put(240, 106){$\alpha$}
\put(208, 30){$\beta$}
\put(268, 32){$\gamma$}
\put(148, 32){$p$}
\put(170, 40){\circle*{3}}
\put(167, 30){$\xi$}
\put(336, 34){$w$}
\put(248, 144){$\lambda$}
\put(178, 130){$\mu$}
\put(100, 40){\vector(3,2){170}}
\put(190, 100){\vector(-2,3){20}}
\put(120, 100){\line(1,0){250}}
\put(130,10){\line(2,3){60}}
\put(350,10){\line(-2, 3){60}}
\thicklines
\put(150, 40){\line(1, 0){180}}
\put(150, 40){\line(3,2){90}}
\put(330, 40){\line(-3,2){90}}
\end{picture}

In order to rule out this subcase we will use arguments of the type in the proofs
of Theorem 6.4 and Theorem 6.12 of \cite{DW6}. First note that $\{\bar{c}, \bar{a},
\bar{p}, \bar{w} \} \subset {\mathcal C}$. Consider the wedge bounded by
the rays $\bar{c}\bar{p}$ and $\bar{c}\bar{a}$. All elements of $\mathcal C$ must
lie in this wedge. Also, any element in the interior of the wedge has positive inner
product with $\bar{c}$. The only elements of $\mathcal C$ lying in the ray $\bar{c}\bar{p}$
are $\bar{c}$ and $\bar{p}$ and the only elements of $\mathcal C$ lying in the ray
$\bar{c}\bar{a}$ are $\bar{c}$ and $\bar{a}$. Similar statements can be made for
the wedge bounded by the rays $\bar{a}\bar{c}$ and $\bar{a}\bar{w}$.

We claim that no element of $\mathcal C$ can lie below $\bar{p}\bar{w}$ in the
intersection of the two wedges. If there is such an element $\bar{b}$, then it must
be null and we must have $J(\bar{b}, \bar{p})=0=J(\bar{b}, \bar{w})$. But the
last two conditions give a contradiction to the nullity of $\bar{b}$. Hence besides
$\bar{c}$ and $\bar{a}$, all elements of $\mathcal C$ lie in the triangle (T10).
It follows that in order to account for the term corresponding to $2\bar{\beta}$
in the superpotential equation, we need to have $\bar{c}^{(1)}, \bar{c}^{(2)}$
in the segment $\bar{p}\bar{w}$ such that $2\bar{\beta}=\bar{c}^{(1)}+ \bar{c}^{(2)}.$
The following arguments shows that this is impossible.

Let us deal with the $d=(3,3,3)$ case first. In this case $c=x$ and $a=y$ and
$\alpha a \beta p$ is a parallelogram. If $\bar{p} \neq \bar{\xi} \in {\mathcal C}$
lies in $\bar{p}\bar{\beta}$, then as long as $\bar{\xi} \neq \bar{\beta},$ the point
$\frac{1}{2}(\bar{c}+\bar{\xi})$ must lie above $\bar{p}\bar{\alpha}$ and so
$\bar{c}+\bar{\xi}$ cannot be cancelled by another pair in $\mathcal C$.
(The easiest way to see this is to introduce orthogonal coordinates
$\mu, \lambda$ as shown.) Therefore,
$\bar{p}\bar{\beta}$ has no interior points in $\mathcal C,$ and either
$\bar{\beta} \in {\mathcal C}$ or $2\bar{\beta}= \bar{p} + \bar{\gamma}$ with
$\bar{\gamma} \in {\mathcal C}$. If $\bar{\beta} \in {\mathcal C}$ then since
$J(\bar{p}, \bar{\beta})= \frac{1}{3} \neq 0$, the term corresponding to
$\bar{p}+\bar{\beta}$ in the superpotential equation is unaccounted for.
But then $\bar{\gamma} \in {\mathcal C},$ which cannot hold by an analogous
argument using the segment $\bar{w}\bar{\gamma}$.

If $d=(9,4,4)$, the difference is that $\frac{1}{2}(\bar{c}+\bar{\xi})$ always
lies above $\bar{p}\bar{\alpha}$ (even when $\bar{\xi}=\bar{\beta}$). So
$\bar{\beta} \notin {\mathcal C}$ and the only way to represent
$2\bar{\beta}$ as $\bar{c}^{(1)}+ \bar{c}^{(2)}$ is $\bar{p}+\bar{\gamma}$ with
$\bar{\gamma} \in {\mathcal C}$. By symmetry, however, we also have
$\bar{\gamma} \notin {\mathcal C}$, a contradiction. So (T10) has been eliminated.

\smallskip

\noindent{\bf Triangle (T12):} We consider the subcase where $\bar{c}+\bar{a}=2\bar{\alpha}$
with $J(\bar{a}, \bar{p})=0=J(\bar{c}, \bar{s})$. These conditions imply that
$c_2 = -\frac{d_2}{2}(1+\frac{1}{d_1})$ and $c_3 =\frac{d_2d_3}{2d_2+d_3}$.
On the other hand, the null condition for $\bar{a}$ implies that $c_1 =-1$,
so $c_2 +c_3 =0$. It follows that $d_1 =\frac{d_3+2d_2}{d_3 -2d_2}$.
The null condition for $\bar{c}$ is $1= \frac{1}{d_1} + \frac{d_2}{4}(1+\frac{1}{d_1})^2
+ d_3(\frac{d_2}{2d_2+d_3})^2$, which implies that $d_2=2$ or $3$.
Substituting the expression for $d_1$ into this null condition and simplifying
gives $ d_3^2 +d_2d_3 -4d_3 -8d_2 =0$, which has no integral solutions with
$d_2=2$ or $3$.

\smallskip

\noindent{\bf Rectangle (Q3):} Notice that interchanging the last two coordinates is
a symmetry of the rectangle. We consider the subcase where $\bar{c}+\bar{a} = 2\bar{x}$
with $J(\bar{c}, \bar{v})=0=J(\bar{a}, \bar{y})$. These conditions imply that
$$\frac{c_1}{d_1}=\frac{1}{2}\left(1-\frac{3}{d_1}+\frac{1}{d_2}+\frac{1}{d_3}\right), \,\,\,
  \frac{c_2}{d_2}=\frac{1}{4}-\frac{5}{4d_1}+\frac{7}{4d_2} + \frac{7}{4d_3}, \,\,\,
  \frac{c_3}{d_3}=-\frac{1}{4}-\frac{3}{4d_1}+\frac{1}{4d_2} +\frac{1}{4d_3}.$$
The null condition for $\bar{c}$ becomes
\begin{eqnarray*}
0&=&d_1^2 d_2(d_2-2)d_3^3 +d_1^2 d_2^2 d_3^2 (3d_1-28) + d_1 d_2^3 d_3^2 (d_1 -10) +
     14d_1 d_2^3 d_3 (d_1-5) \\
  & &   + 6d_1 d_2 d_3^3(d_2 -1) +d_1 d_2^2 d_3^2(d_1^2 -40) + \,\,
     \mbox{{\rm positive quantity}}.
\end{eqnarray*}
Hence $d_1 \leq 9$. For each of these values of $d_1$, we look for further restrictions on
$d_2$ or $d_3$ imposed by the null condition. For example, if $d_1=6$, the null
condition can be rewritten as
$$0=54d_2 d_3^3(d_2-2)+d_2^2 d_3^2(27d_3 -384) + \,\mbox{\rm positive quantity}.$$
It follows that $d_3 \leq 14,$ and for each of these values of $d_3$, the
resulting polynomial equation in $d_2$ has no positive integral roots. By similar
arguments, all the remaining cases can be eliminated.

\smallskip

\noindent{\bf Trapezoid (Q8):} We consider the subcase where $\bar{a}+\bar{c}
=2\bar{z}$ with $J(\bar{c}, \bar{y})=0=J(\bar{a}, \bar{p}).$ Notice that
we have $\frac{c_2}{d_2} =-\frac{1}{2}(1+\frac{1}{d_1}+\frac{1}{d_2}+\frac{1}{d_3})$
and the null condition for $\bar{c}$ is $1= d_1 (\frac{c_1}{d_1})^2 +
d_2 (\frac{c_2}{d_2})^2 + d_3 (\frac{c_3}{d_3})^2.$ Now
$|\frac{c_2}{d_2}| > \frac{1}{2},$ so that $d_2 < 4.$ Also, $d_2 \neq 1$,
otherwise $|\frac{c_2}{d_2}| > 1$ and the null condition would be violated.
Now for $d_2=2,$ the null condition for $\bar{c}$ can be written as
$$ 0 = d_1^3 (3d_3-14)^2 +10d_1(d_1-2)d_3^3 + d_1^2 d_3^2 (15d_3 -160)
     + \, \mbox{\rm positive quantity}.$$
So $d_3 \leq 10,$ and for each of these values, the null condition becomes
a polynomial in $d_1$ which does not have any integral roots. Likewise, for $d_2=3,$
the null condition for $\bar{c}$ becomes
$$0=d_1^3(4d_3^2-84d_3 +441)+18d_1(d_1-2)d_3^2 + d_1^2 d_3^2(18d_3 -264)
    + \, \mbox{\rm positive quantity}.$$
It follows that $d_3 \leq 14,$ and for each of these values, the
resulting null condition has no positive integer roots.

\smallskip

\noindent{\bf Quadrilateral (Q9):} This is (T10) with the point $\bar{v}$
added (below $\bar{p}\bar{w}$ in the diagram above). We are able to rule
out the subcases $\bar{c}+\bar{a}= 2\bar{v}, 2\bar{p}$ and $2\bar{w}$ first.
For the subcase $\bar{c}+\bar{a} =2 \bar{\alpha},$ we are left with the two
special configurations $d=(3,3,3)$ and $d=(9, 4, 4)$ of (T10). To rule these
out, recall that $\bar{w} \in {\mathcal C}$, and so if $\bar{v} \in {\mathcal C}$
as well, then $\bar{v}\bar{w}$ would be an edge with no interior points
in $\frac{1}{2}(d+{\mathcal W})$. By Theorem 3.7 in \cite{DW6}, we would have
$J(\bar{v}, \bar{w})=\frac{1}{4}(1 + \frac{2}{d_3})  > 0,$ a contradiction.

As $\bar{v} \notin {\mathcal C},$ we next consider the term corresponding
to $d+v$ in the superpotential equation. Observe that besides $\bar{a}, \bar{c}$
there can be no further null element of $\mathcal C.$ This follows because
if $\bar{b}$ were such an element, $\Delta^{\bar{b}}$ would have to have
a (1A) and a (1B) endpoint and we have already ruled out all situations involving
the other vertices of (Q9). We cannot have $2\bar{v}=d+v= \bar{c}^{(1)}+\bar{c}^{(2)}$
for $\bar{c}^{(i)} \in \frac{1}{2}(d+{\mathcal W})$ since $\bar{v}$ is a vertex.
We also cannot have one $\bar{c}^{(i)} \in \frac{1}{2}(d+{\mathcal W})$
(and the other equal to $\bar{c}$ or $\bar{a}$) since $c_1=a_1=-1$ while $c_1^{(i)} \leq 1$.
So the superpotential equation has no solution.

\smallskip

\noindent{\bf Pentagon (P1):} Since (P1) contains the triangle (T10), we
again must save the subcase $\bar{c}+\bar{a}=2\bar{\alpha}$ until the last.
The two special configurations $d=(3,3,3)$ and $d=(9,4,4)$ can then be
eliminated as for (Q9).

\smallskip

\noindent{\bf Pentagon (P4):} We first consider the subcase where
$\bar{c}+\bar{a} =2\bar{z}$ with $J(\bar{a}, \bar{x})=0=J(\bar{c}, \bar{s}).$
These conditions imply that $\frac{c_3}{d_3} =\frac{1}{2}(1+\frac{1}{d_1}
+\frac{1}{d_2} - \frac{3}{d_3}).$ If $d_3 > 8,$ then $1-\frac{3}{d_3} \geq
\frac{2}{3}$ and so $\frac{c_3}{d_3} > \frac{1}{3}.$ But then
$d_3 (\frac{c_3}{d_3})^2 > 1,$ which contradicts the nullity of $\bar{c}.$
Now for each value of $d_3 \leq 8$, we examine the null condition for $\bar{c}.$
In each case we obtain further upper bounds on either $d_1$ or $d_2$.
For example, if $d_3=2$, the nullity of $\bar{c}$ can be written as
$$ 0=4d_1^3(d_2^2 -12 d_2 + 36) + d_1^2 d_2^2 (6d_2 -98) + 3d_1 d_2^3(d_1 -12)
+ \,\, \mbox{\rm positive quantity}, $$
and we conclude that either $d_1 \leq 11$ or $2 \leq d_2 \leq 16$. We can then
use MAPLE to find the positive integral roots of the polynomials obtained by
further specializing to these values. In this way this subcase can be eliminated.

Having ruled out all other subcases, we arrive at the situation $\bar{c}
+\bar{a} = 2\bar{z}$ with $J(\bar{a}, \bar{p})=0=J(\bar{c}, \bar{s}).$
Now we know that $\bar{c}, \bar{a}$ are the only null elements in
$\mathcal C$ and $\bar{p}, \bar{s} \in {\mathcal C}.$ By similar arguments
as before, to account for the term $d+y$ in the scalar curvature function,
we need $\bar{y} \in {\mathcal C}.$ But Proposition 3.7 in \cite{DW6}
gives $0=1+\frac{1}{d_1}-\frac{2}{d_2},$ so $d_2=1.$ This contradicts
the presence of the point $s,$ and completes the analysis of (P4).

\smallskip

\noindent{\bf Pentagon (P5):} Having eliminated all other cases, we
consider the final subcase $\bar{c}+\bar{a} = 2\bar{z}$ with
$J(\bar{a}, \bar{p})=0=J(\bar{c}, \bar{w}).$ These conditions and
$c_1+c_2 +c_3 =-1$ give $0=(d_2 -d_3)(d_2d_3 + 9d_1 +4d_2 +4d_3)$, so
$d_2 =d_3.$ Using this in the null condition of $\bar{c}$ yields
$2 \leq d_2 \leq 6$. One then obtains the unique solution
$d=(27, 3, 3)$ with $c=(1, -\frac{5}{3}, -\frac{1}{3})$ and $a=(1, -\frac{1}{3},
-\frac{5}{3}).$ Note that $J(\bar{p}, \bar{x})=0=J(\bar{w}, \bar{y}).$

The only null elements in $\mathcal C$ are $\bar{c}$ and $\bar{a}.$
Arguing as in the last paragraph of (Q9), we see that
$\bar{x}, \bar{y} \in {\mathcal C}.$ It follows that the only way
to express $d+x$ as the sum of two elements in $\mathcal C$ is
as $\bar{x}+\bar{x}$. But $J(\bar{x}, \bar{x})=\frac{1}{4}(1-\frac{1}{27}-\frac{2}{3}) > 0.$
This forces the coefficient $A_{\bar{x}}$ in the scalar curvature
function (\ref{SCF2}) to be positive. This is a contradiction as
$x$ is of type III.

\smallskip

\noindent{\bf Hexagon (H3):} We can again eliminate all subcases, arriving
at the final situation where $\bar{c}+\bar{a}= 2\bar{z}$ and
$J(\bar{c}, \bar{w})=0=J(\bar{a}, \bar{v}).$ As in the previous cases,
the only null elements of $\mathcal C$ are $\bar{c}$ and $\bar{a}.$
We also have $\bar{v}, \bar{w} \in {\mathcal C}.$ So as before we
conclude that $\bar{y} \in {\mathcal C}$. By Proposition 3.7 in
\cite{DW6}, we have $1=\frac{2}{d_2}+\frac{1}{d_3}.$

Next we claim that $\bar{p} \in {\mathcal C}$ as well. Otherwise,
since $\bar{p}$ is a vertex of (H3), $2\bar{p}$ would have to be the sum
of either $\bar{a}$ or $\bar{c}$ with an element of
${\mathcal C} \cap \frac{1}{2}(d+ {\mathcal W}),$ which is not possible.
With $\bar{p} \in {\mathcal C}$, we can apply Proposition 3.7 in
\cite{DW6} again, this time getting $0=1-\frac{4}{d_3}.$ Substituting
$d_3 =4$ in the last equation of the previous paragraph, we see that
$d_2$ cannot be an integer, which is a contradiction.

\medskip

We have finished the analysis of (1A)(1B) endpoints and hence the proof
of Theorem \ref{smallrthm}.

\section{\bf Application to the Exceptional Aloff-Wallach Spaces}

In this section we will illustrate how the results described in
Section 1 can be applied in a concrete situation by examining
the Aloff-Wallach spaces \cite{AW}. For the two exceptional Aloff-Wallach
spaces, which have multiplicities in their isotropic representations,
there are still unresolved issues regarding the existence
and uniqueness of superpotentials (of scalar curvature type), even though
some information can be obtained from Theorem \ref{smallrthm} by
adding discrete symmetries.

Recall that the Aloff-Wallach spaces are the homogeneous manifolds $N_{k, l} :=
{\rm SU}(3)/{\rm U}(1)_{kl}$ where the subgroup ${\rm U}(1)_{kl}$
is the diagonally embedded circle in ${\rm SU}(3)$
with diagonal entries $(e^{ik\theta}, e^{il\theta}, e^{im\theta})$ where
$k, l, m$ are integers with zero sum. For simplicity we will assume
that $k, l$ are relatively prime, so that $N_{k,l}$ is simply connected.
The spaces corresponding to permutations of the same three integers are
equivariantly diffeomorphic. We will therefore make statements with the
understanding that equivalent statements apply to the other diffeomorphic spaces.

The two exceptional Aloff-Wallach spaces are $N_{1,-1}$ and $N_{1, 1}$,
at least from the point of view of the space of homogeneous metrics.
For the generic $N_{kl}$, the space of ${\rm SU}(3)$-invariant metrics
has dimension $4$ since the isotropy representation is the sum of $4$ distinct
irreducible real ${\rm U}(1)$ representations (cf \cite{Wg}). Superpotentials for
the cohomogeneity one Ricci-flat system with $N_{kl}$ as principal orbit
were studied in \cite{CGLP1}, \cite{CGLP2} and \cite{KY}. These authors
showed that there is a superpotential of scalar curvature type (with no null
vectors in $\mathcal C$) and that the first order subsystem associated
to the superpotential corresponds to the ${\rm Spin}(7)$ condition.
Theorem \ref{nonnullthm} shows that there are no further superpotentials
of scalar curvature type without null vectors in $\mathcal C$ while
Theorem \ref{nullthm} shows that there are also no superpotentials of
scalar curvature type with null vectors in $\mathcal C$.

In order to discuss the two exceptional cases, we need to fix some
notation. Let $G={\rm SU}(3), K={\rm U}(1)_{kl}$ and $T$ denote
the set of all diagonal matrices ${\rm diag}(e^{i\theta_1}, e^{i\theta_2}, e^{i\theta_3})$
in $G$. Using the bi-invariant metric $-\tr(XY)$ on the Lie algebra
$\g=\su(3)$, we may decompose the isotropy representation as
\begin{equation} \label{isorep}
 \g/\kf \approx \p = \p_0 \oplus \p_1 \oplus \p_2 \oplus \p_3,
\end{equation}
where $\tf = \kf \oplus \p_0$, $\p_1$ corresponds to the root space
of $\theta_1 -\theta_2$, $\p_2$ corresponds to the root space
of $\theta_1-\theta_3$ and $\p_3$ corresponds to the root space
of $\theta_2-\theta_3$ (cf \cite{Wg}).

\noindent{\bf A. The $N_{1,-1}$ case.}

The summands $\p_2$ and $\p_3$ become equivalent real (irreducible)
representations of $K$. Indeed, upon complexifying $\p_2,$ for example, we get
$\varphi \oplus \varphi^*$ where $\varphi$ is the
standard one-dimensional representation of the circle. It follows
that the space of ${\rm SU}(3)$-invariant metrics has dimension
$6$. The normaliser of $K$ in $G$ is $T,$ and $T/K$ acts non-trivially on
the space of invariant metrics, reducing the effective number
of parameters to $5$, including homothety. By a computation, the
$4$-parameter family of invariant metrics diagonal with respect to
the decomposition (\ref{isorep}) can be shown also to have Ricci tensor
diagonal with respect to this decomposition. Therefore, the superpotential
found in \cite{CGLP1}, \cite{CGLP2} and \cite{KY} using a generic choice
of $k,l$ specialises to a superpotential when $k=1, l=-1$. However, the polytopes
${\rm conv}(\mathcal W)$ and \wt now each have one fewer vertex.

Since the no-multiplicities assumption in Theorems \ref{nonnullthm} and
\ref{nullthm}  is no longer satisfied, it is open whether the known
superpotential is unique among superpotentials of scalar curvature type
(modulo an overall negative sign and an additive constant).

We may, however, add a finite group of isometries to $G$ and $K$ to
eliminate the multiplicities in the isotropy representation. Let
$\hat{G} = G \times (\Z/2)$ and $\hat{K} = K \ltimes \Gamma$
where $\Gamma$ is the diagonally embedded $\Z/2$ such that its image
in ${\rm SU}(3)$ is an order two element of $N_G(T)$ inducing the interchange
of $\theta_1$ and $\theta_2$. The isotropy representation of $\hat{G}/\hat{K}$
now has three irreducible summands given by $\p_0, \p_1$ and $\tilde{\p_2} =
\p_2 \oplus \p_3$. The $\hat{G}$-invariant metrics consist of the $3$-dimensional
subfamily of the $G$-invariant metrics diagonal with respect to (\ref{isorep})
where the scalings along the summands $\p_2$ and $\p_3$ are equal.
Thus $r=3$ and $d=(1, 2, 4)$. The set $\mathcal W$ of weights for the scalar
curvature function is $\{\gamma, \beta, p, v\} = \{(0, -1, 0),\, (0, 0, -1), \, (0, 1, -2),
(1, 0, -2)\}$. The polygon ${\rm conv}(\mathcal W)$ is precisely the triangle
(T4) in Section 3.

Now Theorem \ref{smallrthm}, which allows for disconnected transitive groups,
 implies that the Ricci-flat system for
the $3$-parameter family of invariant metrics on $\hat{G}/\hat{K}$ has no
superpotentials of scalar curvature type with null vectors in $\mathcal C$.
Notice also that $J(\bar{v}, \bar{\gamma}) = \frac{1}{4} \neq 0$. Hence by
Theorem 3.5 in \cite{DW4} there are also no superpotentials of scalar curvature
type without null vectors in $\mathcal C$. In particular this example shows
that the property of having a superpotential of scalar curvature type is
not preserved when we restrict ourselves to a subfamily of invariant metrics
on the principal orbit.

\noindent{\bf B. The $N_{1, 1}$ case.}

In this case, $\p_0 \oplus \p_1$ becomes a $3$-dimensional trivial representation
in the isotropy representation, while $\p_2$ and $\p_3$ become equivalent
$2$-dimensional irreducible summands. The space of ${\rm SU}(3)$-invariant metrics
on $N_{1, 1}$ is now $10$-dimensional. The normalizer of $K$ is
${\rm S}({\rm U}(2){\rm U}(1))$ so that $N(K)/K \approx {\rm SO}(3).$ This
latter group acts on $S^2(\p_0 \oplus \p_1)^K \approx S^2(\R^3)$ by conjugation
(i.e., by the usual action of ${\rm SO}(3)$ on the space of symmetric $3 \times 3$
matrices), and on $S^2(\p_2 \oplus \p_3)^K$ as $\R^3 \oplus \I$, where $\I$ denotes the
trivial representation and $\R^3$ the vector representation.

Again, a superpotential was found in \cite{CGLP1}, \cite{CGLP2}, \cite{KY} for the
Ricci-flat system associated to the $4$-dimensional family of ${\rm SU}(3)$-invariant
metrics diagonal with respect to the decomposition (\ref{isorep}). In fact, a direct
computation shows that these metrics have Ricci tensors diagonal with respect to
(\ref{isorep}). This time the polytopes ${\rm conv}(\mathcal W)$ and \wt are the
same as those for generic choices of $k, l$. The uniqueness of the known superpotential
(modulo an additive constant and an overall minus sign) is again open.

As before, we may eliminate the multiplicities in the isotropy representation
by adding to $G$ a dihedral group $\Gamma$ of symmetries as described in Remark 2.4
in \cite{DW6}. One easily checks that the order two element of the dihedral
group acts as $-1$ on $\p_0$, preserves $\p_1$, and interchanges $\p_2$ and $\p_3$.
The order $3$ element acts trivially on $\p_0$ and acts by rotations on each $\p_i, i> 0$.
Hence the isotropy representation of $\hat{G}/\hat{K}=(G \times \Gamma)/(K \ltimes \Delta\Gamma)$
becomes again the sum of three irreducible summands $\p_0, \p_1$
and $\tilde{\p_2}:=\p_2 \oplus \p_3$. We have $r=3$ with $d=(1, 2, 4)$.
The $G \times \Gamma$ invariant metrics consist of the $3$-parameter family of
diagonal $G$-invariant metrics with the parameters associated to $\p_2$ and $\p_3$
equal.  By specialising the scalar curvature function of the diagonal metrics
we see that $\mathcal W$ is given by $\{s, v, p, \beta, \gamma\} =
\{(1, -2, 0), (1, 0, -2), (0, 1, -2), (0, 0, -1), (0, -1, 0)\}$. Hence the polygon
${\rm conv}(\mathcal W)$ is the parallelogram (Q1) in Section 3.

By Theorem \ref{smallrthm}, the Ricci-flat system associated to
$(G \times \Gamma)/(K \ltimes \Delta\Gamma)$ does not have a superpotential
of scalar curvature type with a null vector in $\mathcal C$. Notice that
$\mathcal W$ is the same as that of case (4) of Theorem \ref{nonnullthm},
but since the groups are not connected, the homogeneous space does
not appear in that classification theorem. In fact one can check that there
is a superpotential with no null elements in $\mathcal C$ given by
$$ u =\sqrt{2}\left(-e^{q_0+2q_2}+e^{q_0+q_1+q_2} +2e^{\frac{q_0}{2}+\frac{3q_1}{2}+q_2}
+4e^{\frac{q_0}{2}+\frac{q_1}{2} + 2q_2}   \right). $$
Furthermore, the associated first order system becomes (2.8) in \cite{KY}
with $a=b$ and $(\alpha_A, \beta_A, \gamma_A)=(1, 1, -2)$ if we set
$4f^2=e^{q_0}, a^2=e^{q_2},$ and $c^2=e^{q_1}$. Thus the superpotential
of non-null type is associated with ${\rm Spin}(7)$ holonomy, but as
observed in \cite{CGLP2} and \cite{KY} there are no solutions of this
first order system that extend smoothly over a special orbit.

\end{document}